\documentclass[11pt]{article}
\usepackage[utf8]{inputenc}
\usepackage{amsmath}
\usepackage{amsfonts}
\usepackage{color}
\usepackage{amssymb}
\usepackage{mathrsfs}
\usepackage{esint}
\usepackage[colorinlistoftodos]{todonotes}
\usepackage[colorlinks=true]{hyperref}
\hypersetup{
	colorlinks=blue,%
	citecolor=red,%
	filecolor=black,%
	linkcolor=blue,%
	urlcolor=blue
}
\usepackage[margin=1.3in]{geometry}

\usepackage[amsmath,thmmarks,hyperref]{ntheorem}
{
	\theoremstyle{nonumberplain}
	\theoremheaderfont{\bfseries}
	\theorembodyfont{\normalfont}
	\theoremsymbol{\mbox{$\Box$}}
	\newtheorem{pf}{Proof.}
}

\numberwithin{equation}{section}
\def\R{\mathbb{R}}
\def\B{\mathbb{B}}
\def\sB{\mathscr{B}^3_3}
\def\Bn{\mathbb{B}^n}
\def\S{\mathbb{S}}

\def\N{\mathbb{N}}

\def \no{\nonumber}
\def\e{\epsilon}
\def\ve{\varepsilon}
\newcommand{\ud}{\mathrm{d}}
\def\pa{\partial}
\newtheorem{thm}{Theorem}[section]

\newtheorem{lem}{Lemma}[section]
\newtheorem{rem}{Remark}[section]

\newtheorem{prop}{Proposition}[section]

\newtheorem*{thm A}{Theorem A}
\newtheorem*{thm B}{Theorem B}
\newtheorem{cor}{Corollary}[section]
\usepackage{upgreek}
\usepackage{graphicx}

\makeatletter

\newdimen\bibspace
\renewenvironment{thebibliography}[1]{%
	\section*{\refname 
		\@mkboth{\MakeUppercase\refname}{\MakeUppercase\refname}}%
	\list{\@biblabel{\@arabic\c@enumiv}}%
	{\settowidth\labelwidth{\@biblabel{#1}}%
		\leftmargin\labelwidth
		\advance\leftmargin\labelsep
		\itemsep\bibspace
		\parsep\z@skip     %
		\@openbib@code
		\usecounter{enumiv}%
		\let\p@enumiv\@empty
		\renewcommand\theenumiv{\@arabic\c@enumiv}}%
	\sloppy\clubpenalty4000\widowpenalty4000%
	\sfcode`\.\@m}
{\def\@noitemerr
	{\@latex@warning{Empty `thebibliography' environment}}%
	\endlist}

\makeatother

\makeatletter

\allowdisplaybreaks

\title{article}
\begin{document}
	\title{\bf  A sharp Sobolev trace inequality of order four \\on three-balls} 
	\author{Xuezhang Chen\thanks{X. Chen: xuezhangchen@nju.edu.cn. Both authors are partially supported by NSFC (No.12271244). }~
		and Shihong Zhang\thanks{S. Zhang: dg21210019@smail.nju.edu.cn.}\\
		{\small $^{\ast}$$^{\dag}$Department of Mathematics \& IMS, Nanjing University, Nanjing 210093, P. R. China }}
	\date{}
	\maketitle

	{\noindent\small{\bf Abstract:}
		We establish a fourth order sharp Sobolev trace inequality on three-balls, and  its equivalence to a third order sharp Sobolev inequality on two-spheres. 		
		
		\medskip 
		
		{{\bf $\mathbf{2020}$ MSC:} 58J32, 46E35 (53A30, 31B30)}
		
		

		\section{Introduction}

		 The geometric inequality on balls and spheres has a long history. We are interested in  the conformally covariant Sobolev (trace) inequality. In a celebrated paper  \cite{Ache-Chang}, Ache and S.-Y. A. Chang established fourth order sharp Sobolev trace inequalities on the unit ball $\B^n$ for  $n\geqslant 4$, which are natural counterparts of inequalities by Lebedev-Milin \cite{Lebedev-Milin} and Beckner \cite{Beckner}.

		For readers' convenience, we restate Ache-Chang sharp Sobolev trace inequalities on $\B^n$ for $n\geqslant 5$.
		\begin{thm A}\label{Thm:Ache-Chang}
			 Let $u\in C^{\infty}(\mathbb{S}^{n-1})$ and $n \geqslant 5$. Then for all $U \in C^\infty(\overline{\Bn})$ satisfying
			 \begin{equation}\label{cond:Neumann_dim_above_three}
			 U=u\quad \mathrm{~~and~~}\quad  \frac{\partial U}{\partial r}=-\frac{n-4}{2}u \qquad \mathrm{~~on~~} \quad \S^{n-1},
			 \end{equation}
			 there holds
			\begin{align}\label{ineq:Ache-Chang_n>4}
			&c_n |\mathbb{S}^{n-1}|^{\frac{3}{n-1}} \left(\int_{\mathbb{S}^{n-1}}|u|^{\frac{2(n-1)}{n-4}}\ud V_{\S^{n-1}}\right)^{\frac{n-4}{n-1}}\no\\
			\leqslant& \int_{\mathbb{B}^n} \left(\Delta U\right)^2 \ud x + 2\int_{\mathbb{S}^{n-1}}|\nabla u|_{\S^{n-1}}^2 \ud V_{\S^{n-1}} +b_n \int_{\mathbb{S}^{n-1}} u^2 \ud V_{\S^{n-1}},
			\end{align}
			where $c_n =n(n-2)(n-4)/4$ and $b_n =n(n-4)/2$. Moreover, equality holds if and only if $U$ is the biharmonic extension of $u_{z_0}(x)=c |1-z_0 \cdot x|^{(4-n)/2}$ on $\S^{n-1}$ and satisfies the Neumann boundary condition, where $c \in \R\backslash\{0\},  z_0 \in \B^n$. 		
		\end{thm A}

		A natural question left in Ache-Chang \cite{Ache-Chang} arises: \emph{Does there exist a sharp Sobolev trace inequality of Ache-Chang type on three-balls?}  The most striking feature is that  it  can reduce to a fractional GJMS equation with a negative critical Sobolev exponent on $\S^2$, which is particularly challenging. The idea of such a reduction traced back to Osgood-Phillips-Sarnak \cite{OPS}, where a derivation of  the Lebedev-Milin inequality from Moser-Trudinger-Onofri inequality was presented. See also  Ache-Chang \cite[p.2739]{Ache-Chang}. Our contribution is to give an affirmative answer to the above question.
				
		\begin{thm}\label{Thm:Sobolev_trace_ineq_three_balls}
				Given  $0<u\in C^{\infty}(\mathbb{S}^{2})$, let $U$ be a smooth extension of $u$ to $\mathbb{B}^3$ satisfying
			\begin{equation}\label{cond:Neumann_dim_three}
			\frac{\partial U}{\partial r} =\frac{1}{2}u\qquad \mathrm{~~on~~}\quad \S^{2},
			\end{equation}
			then
			\begin{align}\label{Ineq:Ache-Chang n=3}
			-\frac{3}{4}|\mathbb{S}^{2}|^{\frac{3}{2}} \left(\int_{\mathbb{S}^{2}}u^{-4}\ud V_{\mathbb{S}^{2}}\right)^{-\frac{1}{2}}
			\leqslant \int_{\mathbb{B}^3} \left(\Delta U\right)^2 \ud x + 2\int_{\mathbb{S}^{2}}|\nabla u|_{\S^{2}}^2 \ud V_{\S^{2}} -\frac{3}{2}\int_{\mathbb{S}^{2}} u^2 \ud V_{\S^{2}},
			\end{align}
			with equality if and only if modulo a positive constant,
			$$U(x)=\sqrt{\frac{|a|^2|x|^2-2a\cdot x+1}{1-|a|^2}}-\frac{1-|x|^2}{4}\sqrt{\frac{1-|a|^2}{|a|^2|x|^2-2a\cdot x+1}}$$
			is biharmonic in $\overline{\B^3}$, where $a \in \B^3$.
		\end{thm}

Theorem \ref{Thm:Sobolev_trace_ineq_three_balls} justifies that  Ache-Chang's  Sobolev trace inequality \eqref{ineq:Ache-Chang_n>4} still holds for $\B^3$.  This combined with Ache-Chang's  inequality draws a complete figure for sharp trace inequalities of order four on balls. As in Ache-Chang \cite{Ache-Chang}, we prefer to use powerful tools in conformal geometry but in a different way, emphasizing the importance of spherical harmonics similar to Beckner \cite{Beckner}.

Next we shall involve a fractional GJMS operator $P_3$ on $\S^2$.
 Although this case is not covered by the scattering theory due to Graham-Zworski \cite{Graham-Zworski}, Branson \cite[Theorem 2.8]{Branson} introduced the fractional GJMS operators $P_{2\gamma}$ on $\S^{n-1}$ for $n\geqslant 3$, as intertwining operators from the viewpoint of representation theory, in the most general case for $\gamma \in \mathbb{C}$ with $-\gamma \notin \frac{n-1}{2}+\mathbb{N}$. In particular, as in \cite{Branson} we introduce
	\begin{align}\label{def:P_3_Branson}
P_3=(B-1)B(B+1) \qquad \mathrm{~~with~~}\quad B=\sqrt{-\Delta_{\S^{n-1}}+\frac{(n-2)^2}{4}},
		\end{align}	
which has the conformal covariance property that for a conformal metric $\hat g=e^{2\tau}g_{\S^{n-1}}$, 
\begin{equation}\label{conf_covariance_P_3}
P_3(e^{\frac{n-4}{2}\tau} \varphi)=e^{\frac{n+2}{2}\tau} \hat P_3(\varphi), \quad \forall~ \varphi \in C^\infty(\S^{n-1}).
\end{equation}
		
\begin{thm}\label{thm:Sobolev_ineq}
		The conformal invariant 
		\begin{align}\label{conf_invariant_ineq}
		Y_3^+(\S^2)=\inf_{0<u\in H^{3/2}(\S^2)}(E[u] \cdot \|u^{-1}\|^2_{L^{4}(\S^2)})
		\end{align}
		is achieved by a smooth positive function $u(x)=c |x-a|^{-1}$ on $\S^2$, where $c \in \R_+,  a \in \B^3$, which together with $c=1$ solves
		\begin{align*}
			P_3u=-\frac{3}{8}u^{-5} \qquad \mathrm{~~on~~}\quad \S^2.
		\end{align*}
	\end{thm}	
	
	We would like to point out that the Beckner's inequality on two-spheres is absent.
	
	A closely related topic is the Sobolev inequality associated to the Paneitz operator on three-manifolds, especially  three-spheres. See \cite{Hang-Yang,Xu,Yang-Zhu,Zhang} etc. Due to extra difficulties arising from the fractional GJMS operator $P_3$, some new techniques have to be developed.
		
		A final step to complete the proof of Theorem \ref{Thm:Sobolev_trace_ineq_three_balls} is the transition from the extremal function on $\S^2$ to its biharmonic extension on $\B^3$. Our unified approach, which is of geometric favor, can be also used to determine extremal functions on balls of Ache-Chang's inequalities, which was recently studied by Ndiaye and L. Sun \cite{Ndiaye-Sun} using a different method.
		
		The paper is organized as follows. In Section \ref{Sect2}, we present some preliminary results of conformal boundary operators, and give an elementary proof of the intimate connection between $P_3$ and an extrinsic GJMS operator  $\mathscr{B}_3^3$ in a class of functions on $\Bn$, which is of independent interest. Just for consideration of notations, we postpone the outline of proof of Theorem \ref{Thm:Sobolev_trace_ineq_three_balls} including the equivalence of inequalities \eqref{Ineq:Ache-Chang n=3} and \eqref{conf_invariant_ineq} to Section \ref{Sect3}. The proof of Theorem \ref{Thm:Sobolev_trace_ineq_three_balls}  occupies the remaining sections. A delicate analysis is conducted to unveil a hidden secret between  a constrained $\mathscr{B}_3^3$ on $\B^3$ and $(-\Delta)^{3/2}$  in $\R^2$. Section \ref{Sect5} is devoted to the extremal problem \eqref{conf_invariant_ineq} on two-spheres. In Section \ref{Sect6}, we determine the explicit extremal functions on the unit ball of \eqref{Ineq:Ache-Chang n=3} and  Ache-Chang's inequalities.
		
		\section{Background}\label{Sect2}
		
		To be self-contained, we collect some basic facts about  conformally covariant boundary operators, as these emerged in various literatures. Besides this, a follow-up paper of the same authors is closely related to these conformal boundary operators.
				
		To continue, we set up some notation. For $n \geqslant 3$, we define 
		\begin{align*}
			\mathcal{N}=\left\{U \in C^\infty(\overline{\Bn})\Big|~~ U=u\quad \mathrm{~~and~~}\quad  \frac{\partial U}{\partial r}=-\frac{n-4}{2}u \quad \mathrm{~~on~~} \quad \S^{n-1}\right\}.
			\end{align*}
			Throughout the paper, let $B_{r}(x_0)$ denote a geodesic ball of radius $r$ and center at $x_0$ in space forms: $\S^2, \R^2$ or $\R^3$. Denote by $\R^n_+:=\{z=(z',z_n) \in \R^n\big|~~z' \in \R^{n-1}, z_n>0\}$ the upper half-space.
Denote by  $I: \R^2 \to \S^2 \backslash \{S\}$ the inverse of stereographic projection, where $S$ is the south pole.
		\subsection{Conformally covariant boundary operators}
		 Suppose $(M,g)$ is a smooth Riemannian manifold of dimension $n \geqslant 3$ with boundary $\pa M$ and $\bar g=g|_{T\partial M}$. Let $R_g$ and $\mathrm{Ric}_g$ be the scalar and
Ricci curvatures.  The second fundamental form  is $\pi(X,Y)=\langle \nabla_X \nu_g,Y\rangle$  and its trace-free part is $\mathring{\pi}(X,Y)=\pi(X,Y)-h_g\langle X,Y\rangle$ for $X,Y \in T\pa M$,  where  $\nu_g$ is the outward unit normal on $\pa M$.  Denote by $H=(n-1)h_g$ the mean curvature. The following
conformally covariant operator of order four are discovered by Paneitz \cite{Paneitz}:
\begin{align*}
P_4^g=\Delta_g^2-\delta \left[\left(\frac{n^2-4n+8}{2(n-1)(n-2)} R_g g-\frac{4}{n-2} \mathrm{Ric}_g\right)\ud\right]+{n-4 \over 2}Q_g,
\end{align*}
where $\delta$ is the divergent operator, $\ud$ is the exterior
differential.  Branson emphasized the zeroth order term of Paneitz operator $P_4^g$, the  $Q$-curvature for $n\neq 4$ (cf. Fefferman-Graham \cite{Fefferman-Graham} in critical dimension four) defined by
\begin{align*}
Q_g=-{1 \over 2(n-1)}\Delta_g R_g+{n^3-4n^2+16n-16 \over 8(n-1)^2(n-2)^2} R_g^2-{2 \over (n-2)^2}|\mathrm{Ric}_g|^2.
\end{align*}
Under conformal change of metrics $g_\tau=e^{2\tau}g$, there holds
\begin{equation}\label{conf_invariance_Paneitz}
P_4^{g_\tau}(\varphi)=e^{-\frac{n+4}{2}\tau} P_4^g(e^{\frac{n-4}{2}\tau}\varphi),\qquad \forall~ \varphi \in C^\infty(M).
\end{equation}

On four-manifolds with boundary, Chang-Qing \cite[p.341]{Chang-Qing} introduced a third order conformally covariant boundary operator by
\begin{align*}
(P_3^b)_g u=&-\frac{1}{2}\frac{\pa }{\pa \nu_g} \Delta_g u-\overline \Delta \frac{\pa u}{\pa \nu_g}-\frac{2}{3}H \overline \Delta u+\langle\pi,\overline{\nabla}^2 u\rangle\\
&+\frac{1}{3}\langle \overline \nabla H,\overline \nabla u \rangle-\left(\mathrm{Ric}_g(\nu_g,\nu_g)-\frac{R_g}{6}\right)\frac{\pa u}{\pa \nu_g}
\end{align*} 
and its associated $T_3$-curvature
$$(T_3)_g=\frac{1}{12}\frac{\pa R_g }{\pa \nu_g}+\frac{1}{6}R_g H-\langle R(\nu_g,\cdot,\nu_g,\cdot),\pi\rangle+\frac{1}{9}H^3-\frac{1}{3}\mathrm{tr}_{\bar g}(\pi^3)-\frac{1}{3}\overline \Delta H,$$
where $R(\cdot,\cdot,\cdot,\cdot)$ is the Riemann curvature tensor. Moreover, $P_3^b$ and $T_3$ have conformally covariant property that if $g_\tau=e^{2\tau} g$, then
\begin{itemize}
\item
$(P_3^b)_{g_\tau}=e^{-3\tau}(P_3^b)_g;$
\item
$(P_3^b)_g +(T_3)_g=(T_3)_{g_\tau} e^{3\tau}.$
\end{itemize} 

		Regarding the generalization of the Chang-Qing boundary operator $P_3^b$ to dimension $n$, we prefer to the formulae of  conformally covariant boundary operators introduced  by  J. Case \cite{Case}: For $\Psi \in C^\infty(\overline{M})$,
			\begin{align*}
\mathscr{B}_0^3 \Psi=&\Psi;\\
\mathscr{B}^{3}_1 \Psi=&\frac{\partial \Psi}{\partial \nu_g}+\frac{n-4}{2}h_g \Psi;\\
\mathscr{B}^{3}_2 \Psi=&-\overline \Delta \Psi+\nabla^2 \Psi(\nu_g,\nu_g)+(n-1)h_g\frac{\partial \Psi}{\partial \nu_g}+\frac{n-4}{2}T_2^3\Psi;\\
\mathscr{B}^{3}_3 \Psi=&-\frac{\partial}{\partial \nu_g}\Delta_g \Psi-2\overline{\Delta}\frac{\partial \Psi}{\partial \nu_g}-\frac{n-4}{2}h_g\nabla^2\Psi(\nu_g,\nu_g)+\frac{4}{n-2}\langle\mathring{\pi},\overline \nabla^2\Psi\rangle\\
&-\frac{3n-8}{2}h_g\overline\Delta \Psi-2(n-5)\langle\overline \nabla h_g,\overline\nabla\Psi\rangle+S_2^3\frac{\partial \Psi}{\partial \nu_g}+\frac{n-4}{2}T_3^3\Psi.
\end{align*}
Here $A_{ij}=\frac{1}{n-2}(R_{ij}-\frac{R_g}{2(n-1)}g_{ij})$ is the Schouten tensor, $J={\rm tr}_g(A)$,  and
\begin{align*}
S_2^3=&-\frac{3n^2-13n+16}{4}h_g^2+\frac{n-8}{2}A(\nu_g,\nu_g)+\frac{3n-8}{2}\bar J+\frac{1}{2}|\mathring{\pi}|^2;\\
T_2^3=&\bar J-A(\nu_g,\nu_g)+\frac{n-3}{2}h_g^2;\\
T_3^3=&\frac{\partial J}{\partial \nu_g}-2\overline\Delta h_g-\frac{4}{n-2}\langle\mathring{\pi},\bar A\rangle+\frac{n-4}{2}h_gA(\nu_g,\nu_g)\\
&+\frac{3n-4}{2}h_g\bar J+\frac{n}{2(n-2)}h_g|\mathring{\pi}|^2-\frac{n^2-3n+4}{4}h_g^3.
\end{align*}
Moreover, if we let $g_\tau=e^{2\tau}g$, then 
			\begin{equation}\label{conf_invariance_bdry_operators}
(\mathscr{B}_k^3)_{g_\tau}(\Psi)=e^{-\frac{n+2k-4}{2}\tau}(\mathscr{B}_k^3)_g(e^{\frac{n-4}{2}\tau} \Psi), \qquad k=0,1,2,3.
\end{equation}
The discovery of the conformal boundary operator $\mathscr{B}_3^3$ is due to Branson-Gover \cite{Branson-Gover} in non-critical dimension $n\neq 4$, and  later extended by Grant \cite{Grant} to the critical dimension four together with a local formula for $\mathscr{B}_2^3$, see Juhl \cite{Juhl}, Stafford \cite{Stafford},  Gover-Peterson \cite{Gover-Peterson} and J. Case \cite{Case} for other treatments. J. Case \cite{Case} used a different approach to find \emph{all} conformal boundary operators $\mathscr{B}_k^3$ and established the self-adjointness of these involved boundary operators. Readers are referred to \cite{Case} for details.

\medskip

It is remarkable (cf. \cite[Lemma 6.3]{Case}) that for $n=4$, there hold $\sB=2(P_3^b)_g$ and $T_3^3=2(T_3)_g+E_g$, where $E_g=4\langle W(\nu_g,\cdot,\nu_g,\cdot),\mathring{\pi}\rangle+\frac{8}{3}\mathrm{tr}_{\bar g}(\mathring{\pi}^3)$ with $W$ as the Weyl tensor, has the property that $E_{g_\tau}=e^{-3\tau}E_g$. Due to similar reason, the boundary curvatures $T_k^3$ for dimension $n$ are \emph{in general} not unique.

			On the model space $(\B^{n}, \S^{n-1}, |\ud x|^2)$, for $U\in \mathcal{N}$, the third order  boundary operator $\mathscr{B}^{3}_{3}$  becomes 
		\begin{equation}\label{bdry operator_dim_n}
	\mathscr{B}^{3}_{3}U=-\frac{\partial \Delta U}{\partial r}-\frac{n-4}{2}\frac{\partial^2 U}{\partial r^2}-\frac{n}{2}\Delta_{\S^{n-1}}u+\frac{(n-4)(n^2-3n+4)}{4}u.
	\end{equation}

		\subsection{The equivalence of  $P_3$ and $\mathscr{B}_3^3$ in class $\mathcal{N}$}
		Let $\{\mathcal{Y}_k; k \in \N\}$ be a complete  $L^2(\S^{n-1})$-orthonormal basis consisting of  spherical harmonics of degree $k$ as eigenfunctions for $-\Delta_{\S^{n-1}}$, solving  $-\Delta_{\S^{n-1}}\mathcal{Y}_k=\lambda_k \mathcal{Y}_k$ with $\lambda_k=k(k+n-2)$. Here, in order to simplify our presentation  we use $\mathcal{Y}_k$ to denote an orthonormal basis of the space of spherical harmonics of degree $k$.  Then, for each $k$ we may write
		$\mathcal{Y}_k=Y_k\big|_{\S^{n-1}}$
		for a harmonic homogeneous polynomial $Y_k$ of degree $k$ on $\R^n$. In other words, $Y_k=|x|^k \mathcal{Y}_k$.
		
		The following elementary result is standard, for instance, see Stein \cite[p.276]{Stein}. Whereas for readers' convenience, we include its proof here.
		\begin{lem}\label{lem:elementary}
			 Expand each $f\in L^2(\S^{n-1})$ as	
			\begin{align*}
			f=\sum_{k=0}^{+\infty}a_k\mathcal{Y}_k.
			\end{align*}
						Then, $f\in C^{\infty}(\S^{n-1})$ if and only if $a_k=O(k^{-N})$ for every  $N \in \mathbb{Z}_+$ as  $k \to \infty$.
			
		\end{lem}
		\begin{pf}
			If $f\in C^{\infty}(\S^{n-1})$, then $\forall~ N \in \mathbb{Z}_+$ we have
				\begin{align*}
				\int_{\S^{n-1}}\left(-\Delta_{\S^{n-1}}\right)^{N}f \mathcal{Y}_{k} \ud V_{\S^{n-1}}=	\int_{\S^{n-1}}f \left(-\Delta_{\S^{n-1}}\right)^{N}\mathcal{Y}_{k}\ud V_{\S^{n-1}}=a_k \lambda_k^{N}.
				\end{align*}
				This yields
				\begin{align*}
				|a_k|\leqslant \frac{1}{\lambda_k^N}\left|\int_{\S^{n-1}}\left(-\Delta_{\S^{n-1}}\right)^{N}f \mathcal{Y}_{k}\ud V_{\S^{n-1}}\right|\leqslant \frac{C(N,f)}{k^{2N}}\qquad\mathrm{\,\,for\,\,}\qquad k\gg 1.
				\end{align*}
				
				Conversely, if $a_k=O(k^{-N})$ for every  $N \in \mathbb{Z}_+$ as  $k \to \infty$,
				then we claim that for every multi-index $\alpha$, there exists a positive constant $C_{\alpha}$ such that
				\begin{align}\label{est:higer_derivative_harmonic_polynomial}
				\max_{|x|\leqslant 1}\left|\frac{\partial^{|\alpha|}}{\partial x^{\alpha}}Y_{k}(x)\right|\leqslant C_{\alpha}k^{|\alpha|+(n-1)/2}.
				\end{align}
				Recall that $\|\mathcal{Y}_k\|_{L^2(\S^{n-1})}=1, \forall~ k \in \N$. Then for any $\ve>0$ we have
				\begin{align*}
				\int_{|x|\leqslant 1+\ve}|Y_{k}(x)|^2\ud x=&\int_{0}^{1+\ve}\int_{\partial B_r(0)}|Y_{k}(x)|^2\ud \sigma \ud r\\
				=&\int_{0}^{1+\ve}r^{n-1+2k}\ud r \int_{\S^{n-1}}|\mathcal{Y}_{k}(\theta)|^2\ud V_{\S^{n-1}}(\theta)=\frac{(1+\ve)^{n+2k}}{n+2k}.
				\end{align*}
			Fix an arbitrary $x_0\in \overline{\Bn}$. By local estimates of higher order derivatives  for harmonic functions  we obtain
				\begin{align*}
				\left|\frac{\partial^{|\alpha|} Y_{k}}{\partial x^{\alpha}}(x_0)\right|\leqslant& \frac{C_{\alpha}}{\ve^{|\alpha|+\frac{n}{2}}}\left(\int_{|x-x_0|\leqslant \ve}|Y_{k}(x)|^2\ud x\right)^{1/2}\\
				\leqslant& \frac{C_{\alpha}}{\ve^{|\alpha|+n/2}}\left(\int_{|x|\leqslant 1+\ve}|Y_{k}(x)|^2\ud x\right)^{1/2}
				\leqslant \frac{C_{\alpha}}{\ve^{|\alpha|+n/2}}\frac{(1+\ve)^{k+n/2}}{\sqrt{n+2k}}.\end{align*}
				Take $\ve=1/k$, the above inequality becomes
				$$\left|\frac{\partial^{|\alpha|} Y_{k}}{\partial x^{\alpha}}(x_0)\right|\leqslant C_{\alpha}k^{|\alpha|+(n-1)/2}\left(1+\frac{1}{k}\right)^{k+n/2}\leqslant C_{\alpha}k^{|\alpha|+(n-1)/2}.
$$
Hence, for every fixed $\alpha$, choose $N>|\alpha|+n$ such that  the series
				\begin{align*}
				\sum_{k=0}^{+\infty}a_k\frac{\partial^{|\alpha|}Y_{k}}{\partial x^{\alpha}}
				\end{align*}
				is uniformly  convergent on $\S^{n-1}$, which implies $f\in C^{\infty}(\S^{n-1})$.	
		\end{pf}
		
		\begin{prop}\label{Prop:equiv_operators}
			Let $n \geqslant 3$ and $u \in C^\infty(\S^{n-1})$. 
				If we expand
		$$u(x)=\sum_{k=0}^{+\infty}u_k \mathcal{Y}_k(x),$$ 
		then the unique solution to
		\begin{align}\label{biharmonic_with_bdry_on_ball}
			\begin{cases}
		\displaystyle\Delta^2 U= 0 \qquad &\mathrm{in}\qquad \B^{n},\\
		\displaystyle U= u \qquad &\mathrm{on}\qquad \S^{n-1},\\
		\displaystyle \frac{\partial U}{\partial r}=-\frac{n-4}{2}u\qquad &\mathrm{on}\qquad \S^{n-1},
		\end{cases}
		\end{align}
		can be expressed as
		$$\sum_{k=0}^{+\infty}u_kY_k(x)\left[1+\left(\frac{k}{2}+\frac{n-4}{4}\right)(1-|x|^2)\right]\in C^{\infty}(\overline{\Bn}).$$
		\end{prop}
		
		\begin{pf}
				We shall solve PDE \eqref{biharmonic_with_bdry_on_ball} by separation of variables. For each $k$, we seek the unique solution of the form
				$$U_k(x)=Y_k(x)(a|x|^2+b) \qquad \mathrm{~~for~~}\quad a, b \in \R$$
				to
			\begin{align*}
			\begin{cases}
			\displaystyle \Delta^2 U_k=0 &\mathrm{~~in~~}\qquad \B^n,\\
			\displaystyle U_k=\mathcal{Y}_k &\mathrm{~~on~~}\qquad \S^{n-1},\\
			\displaystyle \frac{\partial  U_k}{\partial r}=-\frac{n-4}{2}\mathcal{Y}_k \qquad\qquad\qquad&\mathrm{~~on~~}\qquad \S^{n-1}.
			\end{cases}
			\end{align*}
			
			Using $x\cdot\nabla Y_k(x)=kY_k(x)$ we have
			\begin{align*}
				\Delta^2 U_k=2a(n+2k)\Delta Y_k(x)=0.
			\end{align*}
			Substituting $U_k$ into  these two boundary conditions to show
			\begin{align*}
				a+b=1\qquad \mathrm{and}\qquad (k+2)a+kb=-\frac{n-4}{2}.
			\end{align*}
			Eventually we obtain 
		\begin{align*}
		U_k(x)=Y_k(x)\left[1+\left(\frac{k}{2}+\frac{n-4}{4}\right)(1-|x|^2)\right].
		\end{align*}
		
		Next, we  claim that 
		\begin{align*}
			U(x):=\sum_{k=0}^{+\infty}u_kY_k(x)\left[1+\left(\frac{k}{2}+\frac{n-4}{4}\right)(1-|x|^2)\right]\in C^{\infty}(\overline{\Bn})
		\end{align*}
		is a solution to \eqref{biharmonic_with_bdry_on_ball}.
		
		On one hand, it follows from \eqref{est:higer_derivative_harmonic_polynomial} that  $\forall~ m\in \mathbb{N}$, there holds
		\begin{align*}
			\|Y_k\|_{C^{m}(\overline{\Bn})}\leqslant C_m k^{m+(n-1)/2}.
		\end{align*}
		On the other hand, it follows from  Lemma \ref{lem:elementary} that  $ \forall~N\in \mathbb{Z}_+$ we have
		\begin{align*}
			|u_k|=O(k^{-N}) \qquad \mathrm{~~for~~} k \gg 1.
		\end{align*}
		Combining these facts together, for any $x\in \overline{\Bn}$ and sufficiently large $N$, we have
		\begin{align*}
			|\nabla^m U|(x)\leqslant& \sum_{k=0}^{+\infty}C(n)(1+k)\|Y_k\|_{C^m(\overline{\Bn})}|u_k|\\
			\leqslant &\sum_{k=0}^{+\infty}C(n,m)(1+k) k^{m+(n-1)/2}|u_k|<+\infty.
		\end{align*} 
		Hence,  we invoke Ascoli-Arzela theorem to know $U(x)\in C^{\infty}(\overline{\Bn})$.
		\end{pf}
		
		Branson \cite{Branson lec} clarified the relationship between the fractional GJMS operator $P_3$ and its corresponding scattering operator of Graham and Zworski. See also Ache-Chang \cite[Theorem 4.3]{Ache-Chang}. We give an alternative but elementary proof to show that the extrinsic GJMS operator $\mathscr{B}^{3}_{3}$ in class $\mathcal{N}$ agrees with $2 P_3$, which is a special case in \cite[Theorem 1.4]{Case}. 
		\begin{prop}\label{GJMS:extrinsic_intrinsic}
		On the model space $(\B^{n}, \S^{n-1}, |\ud x|^2)$ for $n\geqslant 3$, we have
		 \begin{enumerate}
		\item [(1)] $\mathscr{B}^{3}_{3}U=2 P_3 u$ for all $U\in \mathcal{N}$ satisfying $\Delta^2 U=0$ in $\B^n$ and $u=U\big|_{\S^{n-1}}$.
		\item[(2)] When $n \neq 4$, the fractional $Q$-curvature is $Q_3=\frac{2}{n-4}P_3(1)=\frac{n(n-2)}{4}$.
		\end{enumerate}
		\end{prop}
		\begin{pf}
		For all $U\in \mathcal{N}$, the third order boundary conformally covariant operator becomes
			\begin{align*}
	\mathscr{B}^{3}_{3}U=-\frac{\partial \Delta U}{\partial r}-\frac{n-4}{2}\frac{\partial^2 U}{\partial r^2}-\frac{n}{2}\Delta_{\S^{n-1}}u+\frac{(n-4)(n^2-3n+4)}{4}u
	\end{align*}
		 and the fractional GJMS operator $P_3$ is given in \eqref{def:P_3_Branson}.

			By Proposition \ref{Prop:equiv_operators}, we know that
			\begin{align*}
			U_k(x)=Y_k(x)\left[1+\big(\frac{k}{2}+\frac{n-4}{4})(1-|x|^2\big)\right]
			\end{align*}
			is the unique smooth solution to
			\begin{align*}
			\begin{cases}
			\displaystyle \Delta^2 U=0 \qquad \qquad \qquad &\mathrm{~~in~~}\qquad \Bn,\\
			\displaystyle U=\mathcal{Y}_k &\mathrm{~~on~~}\qquad \S^{n-1},\\
			\displaystyle \frac{\partial  U}{\partial r}=-\frac{n-4}{2}\mathcal{Y}_k &\mathrm{~~on~~}\qquad \S^{n-1}.
			\end{cases}
			\end{align*}
			Using $x\cdot\nabla Y_k(x)=kY_k(x)$ a direct computation yields
			\begin{align*}
			\Delta U_k(x)=-(k+\frac{n-4}{2})(2k+n)Y_k(x)
			\end{align*}
			and
			\begin{align*}
			\frac{\partial ^2U_k}{\partial r^2}=\left[k(k-1)-(2k+1)(k+\frac{n-4}{2})\right]\mathcal{Y}_k(x) \qquad \mathrm{~~on~~}\quad \S^{n-1}.
			\end{align*}
			Hence,  on $\S^{n-1}$ we arrive at 
			\begin{align*}
			\mathscr{B}_3^3(U_k)=&-\frac{\partial \Delta U_k}{\partial r}-\frac{n}{2}\Delta_{\S^{n-1}}\mathcal{Y}_k-\frac{n-4}{2}\frac{\partial^2U_k}{\partial r^2}+\frac{(n-4)(n^2-3n+4)}{4}\mathcal{Y}_k\\
			=&\left\{(k+\frac{n-4}{2})(2k+n)k+\frac{n}{2}\lambda_k-\frac{n-4}{2}\left[-k^2-(n-2)k-\frac{n-4}{2}\right]\right.\\
			&\quad\left.+\frac{(n-4)(n^2-3n+4)}{4}\right\}\mathcal{Y}_k\\
			=&2\left(k+\frac{n-4}{2}\right)\left(k+\frac{n-2}{2}\right)\left(k+\frac{n}{2}\right)\mathcal{Y}_k\\
			=&2\left(\lambda_k+\frac{(n-2)^2}{4}\right)^{1/2}\left(\lambda_k+\frac{(n-2)^2}{4}-1\right)\mathcal{Y}_k=2 P_3 \mathcal{Y}_k.
			\end{align*}
			This together with Proposition \ref{Prop:equiv_operators} directly implies the first assertion.
			
			For the second assertion,  we take a biharmonic function $U_0=1+\frac{n-4}{4}(1-|x|^2) \in  \mathcal{N}$ for $n \neq 4$ to see that
			\begin{align*}
			\frac{n(n-2)(n-4)}{8}=\frac{1}{2}\mathscr{B}_3^3U_0=P_3(1)=\frac{n-4}{2}Q_3.
			\end{align*}
		\end{pf}

		\section{Fourth order sharp Sobolev trace inequality on three-balls}\label{Sect3}
		
		To make our proof transparent, we would like to explain our strategy first. The complete proof of Theorem \ref{Thm:Sobolev_trace_ineq_three_balls} occupies the rest sections.

			\subsection{Outline of the proof}

		\noindent\textbf{\underline{Strategy of proof of Theorem \ref{Thm:Sobolev_trace_ineq_three_balls}}.}
					We consider a constrained minimization problem
			\begin{align}\label{Variation n=3}
			\inf_{U\in \mathcal{N}}\int_{\B^3}\left(\Delta U\right)^2\ud x,
			\end{align}
			where
			\begin{align*}
			\mathcal{N}=\left\{U\Big|~~ U=u\quad \mathrm{~~and~~}\quad  \frac{\partial U}{\partial r}=\frac{u}{2} \quad \mathrm{~~on~~} \quad \S^2\right\}.
			\end{align*}
			
			A direct method can show that  the minimizer $U_1$ of \eqref{Variation n=3} exits and satisfies 
			\begin{align*}
			\begin{cases}
			\displaystyle \Delta^2 U_1=0 \qquad \qquad&\mathrm{~~in~~}\quad \B^3,\\
			\displaystyle U_1=u &\mathrm{~~on~~}\quad \S^2,\\
			\displaystyle \frac{\partial  U_1}{\partial r}=\frac{u}{2} &\mathrm{~~on~~}\quad \S^2.
			\end{cases}
			\end{align*}
			By Proposition \ref{GJMS:extrinsic_intrinsic}, integrating by parts gives
			\begin{align*}
			&\int_{\mathbb{B}^3} \left(\Delta U_1\right)^2 \ud x + 2\int_{\mathbb{S}^{2}}|\nabla u|_{\S^{2}}^2 \ud V_{\S^{2}} -\frac{3}{2}\int_{\mathbb{S}^{2}} u^2 \ud V_{\S^{2}}\\
			=&\int_{\mathbb{S}^{2}}u \mathscr{B}^{3}_3 U_1\ud V_{\S^{2}}=2\int_{\mathbb{S}^{2}}u P_3 u\ud V_{\S^{2}}.
			\end{align*}
			
		Replacing  $U$ by  $U_1$ on the right hand side of \eqref{Ineq:Ache-Chang n=3}, we are motivated to study 
			\begin{align}\label{Reduced Ineq:Ache-Chang n=3}
		-\frac{3}{8}|\mathbb{S}^{2}|^{\frac{3}{2}} \left(\int_{\mathbb{S}^{2}}u^{-4}\ud V_{\mathbb{S}^{2}}\right)^{-\frac{1}{2}}
		\leqslant  \int_{\mathbb{S}^{2}}u P_3 u\ud V_{\S^{2}}
		\end{align}
		and the associated extremal functions.
		
		 We remind that the discussion above indeed demonstrates the equivalence of  \eqref{Reduced Ineq:Ache-Chang n=3} and \eqref{Ineq:Ache-Chang n=3}.
		
		 Suppose $u$ is a minimizer of the inequality \eqref{Reduced Ineq:Ache-Chang n=3}, and modulo a positive constant, solves
		\begin{equation}\label{PDE:two-spheres}
		P_3 u=-\frac{3}{8} u^{-5} \qquad \mathrm{~~on~~} \S^2.
		\end{equation}
Heuristically, pulling back $I^{\ast}(u^{-4}g_{\S^2})=v^{-4}|\ud y|^2$, that is,
		$$v(y)=u\circ I(y)\sqrt{\frac{1+|y|^2}{2}}, \qquad y \in \R^2$$
		has certain decay at infinity and satisfies
		\begin{equation}\label{PDE:whole-space}
		(-\Delta)^{3/2} v=-\frac{3}{8} v^{-5} \qquad \mathrm{~~in~~} \R^2.
		\end{equation}
		Next we consider the Euler-Lagrange equation of the constrained functional associated to the inequality \eqref{Ineq:Ache-Chang n=3}: For $U\big|_{\S^2}=u$,
		\begin{align*}
		\begin{cases}
		\displaystyle \Delta^2 U=0\qquad\, &\mathrm{~~in~~}\qquad \B^3,\\
		\displaystyle 	\frac{\partial U}{\partial r}=\frac{U}{2}\qquad\qquad\qquad \qquad\,&\mathrm{~~on~~}\qquad \S^2,\\
		\displaystyle \mathscr{B}^{3}_3 U=-\frac{3}{4} U^{-5}\qquad&\mathrm{~~on~~}\qquad \S^2.
		\end{cases}
		\end{align*}
		After some delicate analysis, we are able to show that $v$ also satisfies the following integral equation
		\begin{equation}\label{eqn:Integral}
		v(y)=\frac{3}{16\pi}\int_{\R^2}|y-z|v^{-5}(z)\ud z.
		\end{equation}
		See Theorem \ref{Dim 3 integral equ thm} for a precise statement.
		At this stage, the classification theorem of Yan Yan Li \cite{Li} gives the extremal functions of the above integral equation. A direct but short geometric proof, originating from the first author, Wei and Wu \cite{Chen-Wei-Wu}, can be utilized to determine the extremal functions on three-balls from the ones on two-spheres.
		
		Regarding the inequality part for \eqref{Reduced Ineq:Ache-Chang n=3}, some extra difficulties  arise from the nonlocal operator $P_3$, in comparison of Hang-Yang \cite{Hang-Yang}. Fortunately, a deep insight into the relationship among the equations \eqref{PDE:two-spheres}, \eqref{PDE:whole-space} and \eqref{eqn:Integral} paves the way to a complete proof of Theorem \ref{thm:Sobolev_ineq}, as well as Theorem \ref{Thm:Sobolev_trace_ineq_three_balls}.
		\hfill $\Box$

	\section{A bridge between  $\mathscr{B}_3^3$ and $(-\Delta)^{3/2}$}	
	
	When $n=3$, the situation becomes a little subtle.  For any spherical harmonic $\mathcal{Y}_k(x)$ of degree $k$ on $\S^2$, there holds 
		\begin{align*}
			P_3(\mathcal{Y}_k)=\left(\lambda_k+\frac{1}{4}\right)^{1/2}\left(\lambda_k-\frac{3}{4}\right)\mathcal{Y}_k.
		\end{align*}
		Clearly, $P_3$ on $\S^2$ has a negative eigenvalue $-\frac{3}{8}$.
	
	On $(\B^3,\S^2, |\ud x|^2)$, for $U \in C^\infty(\B^3)$ with $u=U|_{\S^2}$, the conformally covariant  boundary operators become
		\begin{align*}
		\mathscr{B}_1^{3}(U)=\frac{\partial U}{\partial r}-\frac{1}{2}u\qquad \mathrm{~~on~~} \quad \S^{2}
		\end{align*}
		and 
		\begin{align}\label{GJMS:extrinsic_three-ball}
		\mathscr{B}^{3}_3(U)=-\frac{\partial \Delta U}{\partial r}-\frac{3}{2}\Delta_{\S^2}u+\frac{1}{2}\frac{\partial ^2U}{\partial r^2}-u \qquad \mathrm{when}\qquad \mathscr{B}_1^{3}(U)=0.
		\end{align}
		
		The purpose of this section is to build a bridge between the extrinsic GJMS operator $\mathscr{B}_3^3$ in class $\mathcal{N}$ and $(-\Delta)^{3/2}$ in $\R^2$, through the investigation of an extension problem on $\B^3$ and an integral equation in $\R^2$.

		\subsection{An extension problem of $\mathscr{B}_3^3$ in the upper half-space}
		
		In this section, we extend to consider a more general setting: For some qualified candidate $T \in C^\infty(\S^2)$, we can assume the solvability of  positive solutions to
			\begin{align}\label{PDE:prescribed_curvature}
		\begin{cases}
		\displaystyle \Delta^2 U=0\qquad\, &\mathrm{~~in~~}\qquad \B^3,\\
		\displaystyle 	\frac{\partial U}{\partial r}=\frac{U}{2}\qquad\qquad\qquad \qquad\,&\mathrm{~~on~~}\qquad \S^2,\\
		\displaystyle \mathscr{B}^{3}_3 U=2T U^{-5}\qquad&\mathrm{~~on~~}\qquad \S^2.
		\end{cases}
		\end{align}
		As before, we let  $u=U\big|_{\S^2}$. As we shall show,   one among obstructions to the above prescribed curvature problem \eqref{PDE:prescribed_curvature} is  the  Kazdan-Warner type condition.
		\begin{prop}\label{prop:Kazdan-Warner}
		Let $U$ be a smooth positive solution to PDE \eqref{PDE:prescribed_curvature} with $u=U\big|_{\S^2}$, then for any conformal vector field $X$ on $\S^2$, 
		\begin{equation*}
		\int_{\S^2} X(T) u^{-4} \ud V_{\S^2}=0.
		\end{equation*}
		\end{prop}
	\begin{pf}
		We consider a functional
		$$I[U]=\left[\int_{\B^3}\left(\Delta U\right)^2 \ud x + 2\int_{\mathbb{S}^{2}}|\nabla u|_{\S^{2}}^2 \ud V_{\S^{2}} -\frac{3}{2}\int_{\mathbb{S}^{2}} u^2 \ud V_{\S^{2}}\right] \left(\int_{\S^2} T u^{-4} \ud V_{\S^2}\right)^{\frac{1}{2}}$$
		for all $0<U \in \mathcal{N}$. 
		
		Suppose $0<U \in \mathcal{N}$ is a critical point of the  above functional $I$ over $\mathcal{N}$. For any $\Phi \in C^\infty(\overline{\B^3})$, we consider a smooth path in $\mathcal{N}$ through $U$ at $t=0$:
		$$U_t=U+t\Phi-\frac{1}{4}(1-|x|^2) t (\Phi-2x \cdot \nabla \Phi), \qquad \mathrm{~~for~~} |t|\ll 1.$$
		Then, we can apply 
		$$0=\frac{\ud}{\ud t}\Big|_{t=0}I[U_t]$$
		 with little effort to show that modulo a positive constant, $U$ solves PDE \eqref{PDE:prescribed_curvature}.
		
		Let $\varphi_t$ be one-parameter family of conformal transformation on $\S^2$ generated by the conformal vector field $X$. We may construct a new smooth path in $\mathcal{N}$ (still denoted by $U_t$) as follows: Write
		$$(\varphi_t)_\ast(u^{-4}g_{\S^2})=u_t^{-4}g_{\S^2},$$
	let $U_t$ be the unique solution to
	\begin{align*}
	\begin{cases}
		\displaystyle \Delta^2 U_t=0\qquad\, &\mathrm{~~in~~}\qquad \B^3,\\
		\displaystyle  U_t =u_t \qquad&\mathrm{~~on~~}\qquad \S^2,\\
		\displaystyle 	\frac{\partial U_t}{\partial r}=\frac{u_t}{2}\qquad\qquad\qquad \qquad\,&\mathrm{~~on~~}\qquad \S^2.
		\end{cases}
	\end{align*}  
		
		By Proposition \ref{GJMS:extrinsic_intrinsic} and \eqref{GJMS:extrinsic_three-ball}, integrating by parts yields
				\begin{align*}
		&\int_{\B^3}\left(\Delta U_t\right)^2 \ud x + 2\int_{\mathbb{S}^{2}}|\nabla u_t|_{\S^{2}}^2 \ud V_{\S^{2}} -\frac{3}{2}\int_{\mathbb{S}^{2}} u_t^2 \ud V_{\S^{2}}\\
		=&\int_{\B^3}U_t \Delta^2 U_t \ud x+\int_{\S^2} u_t \mathscr{B}^{3}_3 U_t \ud V_{\S^2}\\
		=&2\int_{\S^2} u_t P_3 u_t \ud V_{\S^2}=2E[u_t].
		\end{align*} 	
		
		Since $E[v]$ is conformally invariant, we obtain
		\begin{align*}
		0=\frac{\ud}{\ud t}\Big|_{t=0} I[U_t]=2E[u] \frac{\ud}{\ud t}\Big|_{t=0} \left(\int_{\S^2} T\circ \varphi_t u^{-4} \ud V_{\S^2}\right)^{\frac{1}{2}}.
		\end{align*}
This directly yields
\begin{align*}
		\int_{\S^2} X(T) u^{-4} \ud V_{\S^2}=0.
		\end{align*}
			    \end{pf}
		\medskip
		
		Let $F: (\R^{3}_{+}, |\ud z|^2) \to (\B^3,|\ud x|^2)$:
		$$x=F(z)=-\mathbf{e}_3+\frac{2(z+\mathbf{e}_3)}{|z+\mathbf{e}_3|^2}$$
	denote a conformal map  with the property that
		\begin{align*}
		F^\ast (|\ud x|^2)=\left(\frac{2}{(1+z_3)^2+|z'|^2}\right)^2 |\ud z|^2
			:=U_0(z)^{-4}|\ud z|^2,\quad z=(z',z_3)\in \R^3_+.
		\end{align*}
		Notice that $I=F|_{\pa \R^3_+}$ is the inverse of stereographic projection from $\S^2 \backslash\{S\}$ to $\R^2$. 
		
		We now let 
		$$V=U_0 U\circ F \qquad \Longleftrightarrow \qquad F^{*}(U^{-4}|\ud x|^2)=V(z)^{-4}|\ud z|^2.$$ 
		Notice that
		$$(\mathscr{B}_1^3)_{|\ud z|^2}V=-\frac{\pa V}{\pa z_3}$$
		and
		\begin{align*}
		(\mathscr{B}_3^3)_{|\ud z|^2}(V)=\frac{\partial \Delta V}{\partial z_3}+2(\frac{\pa^2}{\pa z_1^2}+\frac{\pa^2}{\pa z_2^2})\frac{\partial V}{\partial z_3}.
		\end{align*}

		For the above PDE   \eqref{PDE:prescribed_curvature}, using conformal change formulae \eqref{conf_invariance_bdry_operators} and \eqref{conf_invariance_Paneitz} we are natural to consider 
		\begin{align}\label{biharmonic_eq_half space}
		\begin{cases}
		\displaystyle \Delta^2 V=0\qquad\qquad\qquad\qquad\qquad\qquad\qquad\, &\mathrm{in}\,\,\,\qquad \R_{+}^3,\\
		\displaystyle 	\frac{\partial V}{\partial z_3}=0\qquad\qquad\qquad \qquad\qquad\,\,&\mathrm{on}\qquad\partial\R^3_{+},\\
		\displaystyle \frac{\partial \Delta V}{\partial z_3}=2(T\circ F)V^{-5}\qquad\qquad\,\,\,\,\,\,\,&\mathrm{on}\qquad \partial\R^3_{+},
		\end{cases}
		\end{align}
		under constraints that
		\begin{align*}
		\int_{\R^2}V^{-4}(z',0)\ud z'=\int_{\S^2}u^{-4}\ud V_{\S^2}
		\end{align*}
		and
		\begin{align}\label{decay_at_infinity}
		\lim_{|z|\to +\infty}\frac{V(z)}{|z|}=\frac{u(S)}{\sqrt{2}}>0.
		\end{align}

					\subsection{An integral equation}	
					
					For clarity, instead we define $\R_+^3:=\{(x,t)\big|~~x \in \R^2, t>0\}$ and let
		\begin{align}\label{dim 3 v}
				v(x)=V|_{\pa \R_+^3}=u\circ I(x)\sqrt{\frac{1+|x|^2}{2}}.
		\end{align}

		\begin{thm}\label{Dim 3 integral equ thm}
			Suppose $U$ is a smooth solution of PDE \eqref{PDE:prescribed_curvature} with $u=U|_{\S^2}$, and $v$ is defined as \eqref{dim 3 v}, then $v$ satisfies both the integral equation
			\begin{align}\label{Eqn:integral}
				v(x)=-\frac{1}{2\pi}\int_{\R^2}|x-y|T\circ I(y)v^{-5}(y)\ud y
			\end{align}
			and
			\begin{align}\label{PDE:fractional}
				\left(-\Delta\right)^{3/2}v(x)=T\circ I(x)v^{-5}(x)  \qquad \mathrm{~~in~~} \R^2.
			\end{align}
		\end{thm}

		The PDE  \eqref{biharmonic_eq_half space} together with constraint \eqref{decay_at_infinity} suggests us study the following boundary value problem  in $\R^3_{+}$:
		\begin{align*}
		\begin{cases}
		\displaystyle \Delta^2 u(x,t)=0\qquad\qquad\qquad\qquad\, &\mathrm{in}\qquad \,\,\R_{+}^3,\\
		\displaystyle 	\partial_{t}u(x,0)=0\qquad\qquad \qquad\qquad\,\,&\mathrm{on}\qquad \partial\R^3_{+},\\
		\displaystyle\partial_{t} \Delta u(x,0)=-f(x)\qquad\,\,\,\,\,\,\,&\mathrm{on}\qquad \partial\R^3_{+},
		\end{cases}
		\end{align*}
		with a decay at infinity that
		\begin{align*}
		\lim_{|(x,t)|\to +\infty}\frac{u(x,t)}{\sqrt{t^2+|x|^2}}=c>0.
		\end{align*}
		Here  $f\in C^{\infty}(\R^2)$ satisfies that for some constant $a>3$,
		\begin{align*}
		f(x)=O(|x|^{-a}) \qquad\mathrm{~~as~~}\quad |x|\to \infty.
		\end{align*}

				We now introduce the following singular integral
		\begin{align}\label{dim 3 v formula}
		\hat v(x,t)=\frac{1}{4\pi}\int_{\R^2}\sqrt{t^2+|x-y|^2} f(y)\ud y.
		\end{align}
				Then it is not hard to see that $v$ solves 
		\begin{align}\label{PDE:integral_solution}
		\begin{cases}
		\displaystyle \Delta^2 \hat v(x,t)=0 \qquad&\mathrm{in}\qquad~~ \R_{+}^3,\\
		\displaystyle 	\partial_{t}\hat v(x,0)=0 \qquad &\mathrm{on}\qquad \partial\R^3_{+},\\
		\displaystyle\partial_{t} \Delta \hat v(x, 0)=-f(x) \qquad &\mathrm{on}\qquad \partial\R^3_{+}.
		\end{cases}
		\end{align}
		
		Let $X=(x,t)$ and
		\begin{align*}
		\beta=\frac{1}{4\pi}\int_{\R^2} f(y)\ud y,
		\end{align*}
then
		\begin{align*}
		\hat v(x,t)-\beta|X|=\frac{1}{4\pi}\underbrace{\int_{\R^2} \frac{|x-y|^2-|x|^2}{\sqrt{t^2+|x-y|^2}+\sqrt{t^2+|x|^2}}f(y)\ud y}_{\mathcal{I}}.
		\end{align*}

		\begin{lem}\label{lem:singular integral estimate}
			Let $\hat v$ be defined as \eqref{dim 3 v formula}. Then for  $|X|\gg1$, there exists a positive constant $C$ such that
			\begin{align*}
				|\hat v(x,t)-\beta|X||\leqslant C \int_{\R^2}|y||f(y)| \ud y.
			\end{align*}
			
		\end{lem}
		\begin{pf}
			We decompose  $\R^2=A_1\cup A_2\cup A_3$, where 
			\begin{align*}
			A_1=&\left\{y||y|<\frac{|x|}{2}\right\},\qquad \qquad\qquad  A_2=\left\{y||x-y|<\frac{|x|}{2}\right\},\\
			A_3=&\left\{y||y|\geqslant \frac{|x|}{2},|x-y|\geqslant \frac{|x|}{2}\right\}.
			\end{align*}
			
			Our discussion is divided into two cases.
			
			\medskip
    \noindent\underline{\emph{Case 1.}} $|X|\gg1$ and $t\geqslant |x|/2$.
                       \medskip
    
			We have  $t\gg 1$ and
		\begin{align*}
1	\leqslant 	\frac{|X|^2}{t^2}=\frac{t^2+|x|^2}{t^2}\leqslant 5.
		\end{align*}
		
			In $A_1$,  we know $\frac{|x|}{2}\leqslant |x-y|\leqslant \frac{3|x|}{2}$ and  $t^2+|x|^2\sim t^{2}+|x-y|^2\sim t^2$, hence
			\begin{align}\label{dim 3 lower bound a_1}
		|\mathcal{I}|\leqslant C\left|\int_{A_1}\frac{|y|^2+2|x||y||}{t}|f(y)|\ud y\right|\leqslant C\int_{\R^2}|y||f(y)|\ud y.
			\end{align}
			
			In $A_2$, we have $|x-y|<|x|/2<t$, $|y|\sim |x|$ and  $t^2+|x|^2\sim  t^2$, then 
						\begin{align}\label{dim 3 lower bound a_2}
			|\mathcal{I}|\leqslant	\left|\int_{A_2}\frac{(|x-y|+|x|)|y|}{t}|f(y)|\ud y\right|\leqslant	C\int_{A_2}|y||f(y)|\ud y.
			\end{align}
			
			In $A_3$, we have $|x-y|\leqslant |y|+|x|\leqslant 3|y|$ and $|y|\leqslant |x-y|+|x|\leqslant 3|x-y|$. This directly implies that $|x-y|\sim|y|$. We further decompose $A_3=A_{31}\cup A_{32}$ by
			\begin{align*}
			A_{31}=A_3\cap \{y||y|>t\}, \qquad\qquad A_{32}=A_3\cap \{y||y|\leqslant t\}.
			\end{align*}
			In $A_{31}$, we have $t^2+|x-y|^2\sim |y|^2$ and $t^2+|x|^2\sim t^2$, and then
			\begin{align}\label{dim 3 lower bound b_2}
			|\mathcal{I}|\leqslant C\int_{A_{31}}\frac{(|x-y|+|x|)|y|}{t+|y|}|f(y)|\ud y\leqslant C\int_{A_2}|y||f(y)|\ud y.
			\end{align}
			In $A_{32}$, we have $t^2+|x-y|^2\sim t^2$ and $t^2+|x|^2\sim t^2$, and then
			\begin{align}\label{dim 3 lower bound b_3}
			|\mathcal{I}|\leqslant C\int_{A_{32}}\frac{(|y|+|x|)|y|}{t}|f(y)|\ud y\leqslant C\int_{A_3}|y||f(y)|\ud y.
			\end{align}
			
			Combining \eqref{dim 3 lower bound a_1}, \eqref{dim 3 lower bound a_2}, \eqref{dim 3 lower bound b_2} and \eqref{dim 3 lower bound b_3}, we obtain
			\begin{align*}
				|\hat v(x,t)-\beta|X||\leqslant C  \int_{\R^2}|y||f(y)| \ud y.
			\end{align*}
			
		\medskip
    \noindent\underline{\emph{Case 2.}} $|X|\gg1$ and $t<|x|/2$.
                   \medskip
    
   We have $|x|\gg 1$ and
   \begin{align*}
		1	\leqslant 	\frac{|X|^2}{|x|^2}=\frac{t^2+|x|^2}{|x|^2}\leqslant \frac{5}{4}.
		\end{align*}
		
		In $A_1$,  we have $\frac{|x|}{2}\leqslant |x-y|\leqslant \frac{3|x|}{2}$, then $t^2+|x|^2\sim |x|^2$ and $t^2+|x-y|^2\sim |x|^2$. Thus, 
			\begin{align}\label{dim 3 singular integral estimate a_1}
			|\mathcal{I}|\leqslant C\int_{A_1}\frac{(|x-y|+|x|)|y|}{|x|}|f(y)|dy\leqslant C\int_{A_1}|y||f(y)|\ud y\leqslant C.
			\end{align}
			
			In $A_2$, there holds $|x|/2<|y|<3|x|/2$, which means $|y|\sim |x|$. Further decompose $A_2=A_{21}\cup A_{22}$ by
			\begin{align*}
			A_{21}=A_2\cap \{y||x-y|>t\} \qquad \mathrm{and}\qquad A_{22}=A_2\cap \{y||x-y|\leqslant t\}.
			\end{align*}
			For $A_{21}$,  we have $t^2+|x-y|^2\sim |x-y|^2$ and then
			\begin{align}\label{dim 3 singular integral estimate a_2}
			|\mathcal{I}|\leqslant& \int_{A_{21}}\frac{(|x-y|+|x|)|y|}{|x-y|+|x|}|f(y)|\ud y\leqslant C\int_{A_{21}}|y||f(y)|\ud y.
			\end{align}
			For $A_{22}$, we have $t^2+|x-y|^2\sim t^2$ and then
			\begin{align}\label{dim 3 singular integral estimate a_4}
			|\mathcal{I}|\leqslant &C\int_{A_{22}}\frac{(|x-y|+|x|)|y|}{t+|x|}|f(y)|\ud y\leqslant C\int_{A_{22}}|y| |f(y)|\ud y.
			\end{align}
			
		In $A_3$, we have $|y|/3<|x-y|<3|y|$, $|y|>|x|/2>t$. Thus, $t^2+|x-y|^2\sim |y|^2$, $t^2+|x|^2\sim |x|^2$, there holds
			\begin{align}\label{dim 3 singular integral estimate a_6}
			|\mathcal{I}|\leqslant \int_{A_{3}}\frac{(|x-y|+|x|)|y|}{|x|+|y|}|f(y)|\ud y\leqslant \int_{A_{3}}|y||f(y)|\ud y.
			\end{align}
			
			Combining \eqref{dim 3 singular integral estimate a_1}, \eqref{dim 3 singular integral estimate a_2},  \eqref{dim 3 singular integral estimate a_4} with \eqref{dim 3 singular integral estimate a_6},  we conclude that
			\begin{align*}
			|v(x,t)-\beta|X||\leqslant C  \int_{\R^2}|y||f(y)| \ud y.
			\end{align*}

		\end{pf}

		\begin{thm}\label{dim 3 Integral equ thm}
			Assume $f \in C^\infty(\R^2)$ satisfies $f(x)=O(|x|^{-a})$ as $|x| \to \infty$ for some constant $a>3$. Let $u$ be a smooth solution of 
			\begin{align*}
			\begin{cases}
			\displaystyle \Delta^2 u(x,t)=0\qquad\qquad\qquad\qquad\, &\mathrm{in}\qquad \R_{+}^3,\\
			\displaystyle 	\partial_{t}u(x,0)=0\qquad\qquad \qquad\qquad\,\,&\mathrm{on}\quad\,\,\,\, \partial\R^3_{+},\\
			\displaystyle\partial_{t} \Delta u(x,0)=-f(x)\qquad\,\,\,\,\,\,\,&\mathrm{on}\quad\,\,\,\, \partial\R^3_{+},
			\end{cases}
			\end{align*}
		under the constraint that
			\begin{align}\label{Decay_rate_at_infinity}
		\lim_{|(x,t)|\to +\infty}\frac{u(x,t)}{\sqrt{t^2+|x|^2}}=c>0.
			\end{align}
			Then there exists some constant $C$ such that
			\begin{align*}
			u(x,t)=\frac{1}{4\pi}\int_{\R^2}\sqrt{t^2+|x-y|^2} f(y)\ud y+C.
			\end{align*}
			\end{thm}
		\begin{pf}
			By \eqref{PDE:integral_solution} and Lemma \ref{lem:singular integral estimate}, we know that
			  \begin{align*}
			w(x,t)=u(x,t)-\hat v(x,t)
			\end{align*}
			solves the following PDE
			\begin{align*}
			\begin{cases}
			\displaystyle \Delta^2 w(x,t)=0\qquad &\mathrm{in}\qquad~~ \R_{+}^3,\\
			\displaystyle 	\partial_{t}w(x,0)=0\qquad &\mathrm{on}\qquad\partial\R^3_{+},\\
			\displaystyle\partial_{t} \Delta w(x,0)=0\qquad&\mathrm{on}\qquad \partial\R^3_{+},
			\end{cases}
			\end{align*}
			with decay at infinity that
			\begin{align*}
			\left|w(x,t)\right|=O(|X|) \qquad \mathrm{~~for~~}\quad |X| \gg 1.
			\end{align*}

			We extend $w$ to  $\R^3$ by an even reflection
			\begin{align*}
			\tilde{w}(x,t)=\begin{cases}
			\displaystyle w(x,t)\qquad&\mathrm{for}\qquad t\geqslant 0,\\
			\displaystyle w(x,-t)\qquad&\mathrm{for}\qquad t<0,
			\end{cases}
			\end{align*}
			then $\partial_{t} \tilde{w}(x,0)=0$ and $$\partial_{t}^3w(x,0)=\partial_{t}\Delta w(x,0)-\partial_{t}\Delta_x w(x,0)=0-\Delta_x\partial_t w(x,0)=0.$$
			So, $\tilde{w}$ is biharmonic on $\R^3$. Then the higher derivative estimates for biharmonic functions (for example, see \cite[Proposition 4]{Martinazzi}) implies that 
			\begin{align*}
			\|\nabla^2 \tilde{w}\|_{L^{\infty}(B_R)}\leqslant \frac{C}{R^5}\int_{B_{2R}}|\tilde{w}| \ud X\leqslant \frac{C}{R}  \to 0 \qquad\mathrm{as}\qquad R\to +\infty.
			\end{align*}
			This implies
			$$\tilde w(x,t)=a_1x_1+a_2 x_2+a_3 t+b, \qquad a_i, b \in \R.$$
			On the other hand, by Lemma \ref{lem:singular integral estimate} and \eqref{Decay_rate_at_infinity} we know that the limit $\lim_{|X|\to \infty} \frac{\tilde w(X)}{|X|}$ exists. Hence, it is not hard to see that $a_i=0$ and the desired assertion follows.
		\end{pf}
		
		\subsection{Proof of Theorem \ref{Dim 3 integral equ thm}}
	We now apply Theorem \ref{dim 3 Integral equ thm} to PDE \eqref{biharmonic_eq_half space} together with $v$ as in \eqref{dim 3 v} and obtain
	\begin{align}\label{v_integral sol}
	v(x)=\frac{1}{4\pi}\int_{\R^2}|x-y|f(y)\ud y+C,
	\end{align}
	where
	\begin{align*}
		f(x)=-2T\circ I(x)v^{-5}(x)=O\left(\frac{1}{|x|^5}\right)\qquad \mathrm{for}\qquad |x|\gg 1.
	\end{align*}
	
Recall that
	\begin{align*}
	v(x)=	u\circ I(x)\sqrt{\frac{1+|x|^2}{2}}.
	\end{align*}
	We introduce $I_S: y \in \R^2 \mapsto I(|y|^{-2}y)\in\S^2 \backslash \{N\}$, where $N$ is the north pole. Then near $S$ corresponding to  $|x|=+\infty$, we obtain the following expansion:
\begin{align}\label{Expansion v}
	v(x)=&u\circ I_S\left(\frac{x}{|x|^2}\right)\sqrt{\frac{1+|x|^2}{2}}\nonumber\\
	=&\frac{1}{\sqrt{2}}\left(u\circ I_S(0)+\nabla (u\circ I_S)(0)\cdot\frac{x}{|x|^2}+O\left(\frac{1}{|x|^2}\right)\right)\left(|x|+O\big(\frac{1}{|x|}\big)\right)\nonumber\\
	=&\frac{1}{\sqrt{2}}u(S)|x|+\sum_{i=1}^{2}a_i\theta_i+O\left(\frac{1}{|x|}\right),
\end{align}
where $\theta=|x|^{-1}x$ and 
\begin{align*}
	a_i=\frac{1}{\sqrt{2}}\partial_i(u\circ I_S)(0), \qquad i=1,2.
\end{align*}

We simplify
\begin{align*}
	w(x)=\frac{1}{2\pi}\int_{\R^2}|x-y|f(y)\ud y.
\end{align*}
In the next lemma, we give the expansion of $w$ at infinity.

\begin{lem}\label{Expansion lem w}
	For $|x|\gg 1$ there holds
	\begin{align*}
		w(x)=\alpha|x|+\sum_{i=1}^{2}b_i\theta_i+O\left(\frac{1}{|x|}\right),
	\end{align*}
	where $\alpha=(4\pi)^{-1} \int_{\R^2}f(y)\ud y$.
\end{lem}
\begin{pf}
	Observe that
	\begin{align*}
4\pi(w(x)-\alpha|x|)=&\int_{\R^2}(|x-y|-|x|)f(y)\ud y=\int_{\R^2}\frac{|y|^2-2x\cdot y}{|x|+|x-y|}f(y)\ud y\\
=&-2\int_{\R^2}\frac{x\cdot y}{|x|+|x-y|}f(y)\ud y+O\left(\frac{1}{|x|}\right),
	\end{align*}
	where the last equality follows by $|y|^2|f(y)|\in L^1(\R^2)$. Moreover, the first term can be estimated via
	\begin{align*}
		\int_{\R^2}\frac{x\cdot y}{|x|+|x-y|}f(y)\ud y=\int_{\R^2}\frac{x\cdot y}{2|x|}f(y)\ud y+O\left(\frac{1}{|x|}\right).
	\end{align*}
	This is due to
	\begin{align*}
		&\left|	\int_{\R^2}\left(\frac{1}{|x|+|x-y|}-\frac{1}{2|x|}\right)x\cdot yf(y)\ud y\right|\\
		=&\left|	\int_{\R^2}\frac{|x-y|-|x|}{2(|x|+|x-y|)|x|}(x\cdot y)f(y)\ud y\right|\\
		=&\left|	\int_{\R^2}\frac{|y|^2-2x\cdot y}{2(|x|+|x-y|)^2|x|}(x\cdot y)f(y)\ud y\right|\\
		\leqslant& \int_{\R^2}\frac{|x-y|}{2(|x|+|x-y|)^2} |y|^2 |f(y)|\ud y+\int_{\R^2}\frac{(x \cdot y)^2}{2|x|(|x|+|x-y|)^2} |f(y)|\ud y\\
		\leqslant& \frac{1}{|x|}\int_{\R^2}|y|^2|f(y)|\ud y.
	\end{align*}
	
	Hence we obtain
		\begin{align*}
	w(x)=\alpha|x|+\sum_{i=1}^{2}b_i\theta_i+O\left(\frac{1}{|x|}\right)
	\end{align*}
with
	\begin{align*}
		b_i=-\frac{1}{2\pi}\int_{\R^2}y_if(y)\ud y,\qquad i=1,2.
	\end{align*}
\end{pf}

	\noindent \textbf{Proof of Theorem \ref{Dim 3 integral equ thm}. } By  Lemma \ref{Expansion lem w} and \eqref{Expansion v}, we apply \eqref{v_integral sol}  to know that for $|x|\gg 1$, 
	\begin{align*}
		\frac{1}{\sqrt{2}}u(S)|x|+\sum_{i=1}^{2}a_i\theta_i=\alpha|x|+\sum_{i=1}^{2}b_i\theta_i+C+O\left(\frac{1}{|x|}\right).
	\end{align*}
	Comparing coefficients on both sides, we deduce that $\alpha=\frac{1}{\sqrt{2}}u(S)$, $a_i=b_i, i=1,2$, and also $C=0$. Thus, the desired assertion follows.
	\hfill $\Box$

		\section{Third order sharp Sobolev inequality on two-spheres}\label{Sect5}

		Let $\{\mathcal{Y}_k; k \in \N\}$ be a complete  $L^2(\S^2)$-orthonormal basis consisting of  spherical harmonics of degree $k$ as eigenfunctions for $-\Delta_{\S^2}$, solving  $-\Delta_{\S^2}\mathcal{Y}_k=\lambda_k \mathcal{Y}_k$, where $\lambda_k=k(k+1)$ and
	$$\lambda_0=1<\lambda_1=\lambda_2=\lambda_3=2<\lambda_4 \leq \cdots.$$	
	 The Sobolev space $H^s(\S^2)$ for $0\leqslant s \in \R$ is given by  
		\begin{align*}
			H^s(\S^2)=\left\{u=\sum_{k=0}^{+\infty}u_k\mathcal{Y}_k\Big|~~ \sum_{k=0}^\infty (\lambda_k^s+1) u_k^2<\infty\right\}
		\end{align*}
	coupled with its norm
	$$\|u\|_{H^s(\S^2)}=\sqrt{\sum_{k=0}^\infty (\lambda_k^s+1) u_k^2}.$$

	The following interpolation inequality is needed:
	\begin{equation}\label{ineq:interpolation}
	\|u\|_{H^{1/2}(\S^2)}^2\leqslant C\|u\|_{L^2(\S^2)}\|u\|_{H^{1}(\S^2)}.
	\end{equation}
	This directly follows by
	\begin{align*}
	\|u\|_{H^{1/2}(\S^2)}^2 =\sum_{k=0}^\infty (\lambda_k^{1/2}+1) u_k^2\leqslant& \sqrt{2}(\sum_{k=0}^\infty  u_k^2)^{1/2} (\sum_{k=0}^\infty (\lambda_k+1) u_k^2)^{1/2}\\
	=&\sqrt{2}\|u\|_{L^2(\S^2)}\|u\|_{H^{1}(\S^2)}.
	\end{align*}
	For convenience, we simplify the average of $u$ over $\S^2$ by $\bar u$, and use another equivalent norm of $H^{1/2}(\S^2)$:
		\begin{align*}
		\|u\|_{H^{1/2}(\S^2)}=\left[\int_{\S^2} u(x)\int_{\S^2} \frac{u(x)-u(y)}{|x-y|^3}  \ud V_{\S^2}(y)  \ud V_{\S^2}(x)+\int_{\S^2} u^2 \ud V_{\S^2}\right]^{1/2},
	\end{align*}
	where $|x-y|$ means the distance of $x$ and $y$ in $\R^3$. 
	
		By definition \eqref{def:P_3_Branson} of $P_3$ we have
		\begin{align}\label{P3 expansion}
		P_3u=\sum_{k=0}^{+\infty}\left(\lambda_k+\frac{1}{4}\right)^{1/2}\left(\lambda_k-\frac{3}{4}\right)u_k\mathcal{Y}_k.
		\end{align}
		
		We introduce a quadratic functional on $H^3(\S^2) \times H^3(\S^2)$ by
		\begin{align*}
		E[u,v]=\frac{1}{4\pi}\int_{\S^2}vP_3u   \ud V_{\S^2}.
		\end{align*}	
		In particular, we define the energy for $P_3$ by setting $E[u]=E[u,u]$, and introduce
		\begin{align*}
		Y_3^+(\S^2):=\inf_{0<u\in H^{3/2}(\S^2)}(E[u]\cdot \|u^{-1}\|^2_{L^{4}(\S^2)}).
		\end{align*}
		 Clearly, it follows from \eqref{conf_covariance_P_3}  that  $Y_3^+(\S^2)$ is a conformal invariant.

	A direct consequence of the expansion \eqref{P3 expansion} is that the operator $P_3u+\frac{3}{8}\bar{u}$ is nonnegative. Notice that
	\begin{align*}
		\tilde{E}[u,v]=&\frac{1}{4\pi}\int_{\S^2}  v\left(P_3u+\frac{3}{8}\bar{u}\right) \ud V_{\S^2}=\sum_{k=1}^{+\infty}\left(\lambda_k+\frac{1}{4}\right)^{1/2}\left(\lambda_k-\frac{3}{4}\right)u_kv_k\\
		\leqslant &\left(\sum_{k=1}^{+\infty}\left(\lambda_k+\frac{1}{4}\right)^{1/2}\left(\lambda_k-\frac{3}{4}\right)u^2_k\right)^{1/2}\left(\sum_{k=1}^{+\infty}\left(\lambda_k+\frac{1}{4}\right)^{1/2}\left(\lambda_k-\frac{3}{4}\right)v^2_k\right)^{1/2}\\
		=&\left(E[u]+\frac{3}{8}\bar{u}^2\right)^{1/2}\left(E[v]+\frac{3}{8}\bar{v}^2\right)^{1/2}.
	\end{align*}
	This yields
	\begin{align}\label{Energy inequlity}
		E[u,v]\leqslant \left(E[u]+\frac{3}{8}\bar{u}^2\right)^{1/2}\left(E[v]+\frac{3}{8}\bar{v}^2\right)^{1/2}-\frac{3}{8}\bar{u}\bar{v}.
	\end{align}

	\medskip
	
	Next, our goal is to prove  the following Proposition \ref{Prop: nonnegative_energy}. However, the conformally covariant nonlocal operator $P_3$ brings us some challenges in comparison with \cite{Hang-Yang}.
		
		\begin{lem}\label{V3 lem 1}
			Suppose $u\in H^{3/2}(\S^2)$, $\|u^{-1}\|_{L^4(\S^2)}\leqslant 1$, $\|u\|_{H^{3/2}(\S^2)}\leqslant A$ for some $A \in \R_+$, then there exists a positive constant $c$ such that $u\geqslant c Ae^{-cA^4}$.
		\end{lem}
		\begin{pf}
			By Sobolev embedding theorem we have 
			\begin{align*}
			\|u\|_{C^{1/2}(\S^3)}\leqslant CA.
			\end{align*} 
		Fix an arbitrary $x_0\in \S^2$, then for any $x\in B_{1}(x_0)\subset \S^2$ there holds
			\begin{align*}
			|u(x)|\leqslant |u(x_0)|+Cd_{\S^2}(x_0,x)^{1/2}.
			\end{align*}
			From the assumption,
			\begin{align*}
			1\geqslant& \int_{B_{1}(x_0)}|u(x)|^{-4}\ud V_{\S^2}\geqslant \int_{B_{1}}\left(|u(x_0)|+CAd_{\S^2}(x_0,x)^{1/2}\right)^{-4}\ud V_{\S^2}\\
			\geqslant& C\int_{B_{1}(x_0)}\left(|u(x_0)|^4+CA^{4}d_{\S^2}(x_0,x)^{2}\right)^{-1}\ud V_{\S^2}\\
			\geqslant& \frac{C}{A^4}\log \frac{|u(x_0)|^4+A^4}{|u(x_0)|^4}-C.
			\end{align*}
			This implies the desired assertion.
		\end{pf}

		\begin{lem}\label{Approximate lem}
			Suppose $u\in H^{3/2}(\S^2)$ satisfies $u(S)=0$. Let $\eta_{\ve}\in C^{\infty}(\S^2)$ be a cut-off function  such that $\eta_{\ve}= 1$ in $B_{\ve}(S)$ and $\eta_{\ve}= 0$ in $\S^2\backslash B_{2\ve}(S)$,  then $\eta_{\ve} u\to 0$ in $H^{3/2}(\S^2)$.
		\end{lem}
		\begin{pf}
			We first claim that 
			\begin{align}\label{lim:intermediate_norm}
				\|\eta_{\ve}u\|_{W^{1,4}(\S^2)}\to 0 \qquad \mathrm{~~as~~}\qquad \ve\to 0.
			\end{align}
			To this end, by the Sobolev embedding theorem,  we know $H^{3/2}(\S^2)\hookrightarrow W^{1,4}(\S^2)\hookrightarrow C^{1/2}(\S^2)$, then we use the assumption $u(S)=0$ to estimate
			\begin{align}\label{approximate lem formula 1}
				\left(\int_{\S^2} |\nabla (\eta_{\ve} u)|^4  \ud V_{\S^2}\right)^{1/4}\leqslant& \left(\int_{\S^2} |\nabla \eta_{\ve}|^4u^4\ud V_{\S^2}\right)^{1/4}+\left(\int_{\S^2} |\nabla u|^4\eta_{\ve}^4  \ud V_{\S^2}\right)^{1/4}\nonumber\\
				\leqslant&\ve^{-1}[u]_{C^{1/2}(B_{2\ve}(S))}\left(\int_{B_{2\ve}(S)}d_{\S^2}(S,x)^2\ud V_{\S^2}(x)\right)^{1/4}+\|\nabla u\|_{L^4(B_{2\ve}(S))}\nonumber\\
				\leqslant& C\left([u]_{C^{1/2}(B_{2\ve}(S))}+\|\nabla u\|_{L^4(B_{2\ve}(S))}\right)\nonumber\\
				\leqslant& C\|u\|_{W^{1,4}(B_{2\ve}(S))}\to 0.
			\end{align}
			Clearly, there holds $\|\eta_{\ve} u\|_{L^4(\S^2)}\to 0$ as $\ve\to 0$.  
			
			It remains to estimate $\|\partial_i(\eta_{\ve}u)\|_{H^{1/2}(\S^2)}$. Notice that
			\begin{align*}
				\|\partial_i(\eta_{\ve}u)\|^2_{H^{1/2}(\S^2)}\leqslant C\left(\|\eta_{\ve}\partial_iu\|^2_{H^{1/2}(\S^2)}+\|\partial_i\eta_{\ve}u\|^2_{H^{1/2}(\S^2)}\right).
			\end{align*}

			For the first term, we have
			\begin{align}\label{approximate lem formula 2}
				I_{\e}=&\int_{\S^2}\eta_{\ve}(x)\partial_iu(x)\int_{\S^2} \frac{\eta_{\ve}(x)\partial_iu(x)-\eta_{\ve}(y)\partial_iu(y)}{|x-y|^3}  \ud V_{\S^2}(y)  \ud V_{\S^2}(x)\nonumber\\
				= &\int_{\S^2}\eta_{\ve}(x)\left(\partial_i u(x)\right)^2\int_{\S^2} \frac{\eta_{\ve}(x)-\eta_{\ve}(y)}{|x-y|^3}  \ud V_{\S^2}(y)  \ud V_{\S^2}(x)\nonumber\\
			&+\int_{\S^2}\eta_{\ve}(x)\partial_iu(x) \int_{\S^2} \frac{\eta_{\ve}(y)(\partial_iu(x)-\partial_iu(y))}{|x-y|^3}  \ud V_{\S^2}(y)  \ud V_{\S^2}(x)\nonumber\\
			\leqslant& C\|\nabla u\|^2_{L^4(B_{2\ve}(S))}+\frac{1}{2}\int_{\S^2} \int_{\S^2} \frac{\eta_{\ve}(x)\eta_{\ve}(y)(\partial_iu(x)-\partial_iu(y))^2}{|x-y|^3}  \ud V_{\S^2}(y)  \ud V_{\S^2}(x),
			\end{align}
			where the last inequality follows by 
			\begin{align*}
				\left|\int_{\S^2} \frac{\eta_{\ve}(x)-\eta_{\ve}(y)}{|x-y|^3}  \ud V_{\S^2}(y)\right|\leqslant& C\|\nabla^2\eta_{\ve}\|_{L^{\infty}(\S^2)}\int_{B_{4\ve}(S)}\frac{1}{|x-y|}\ud V_{\S^2}(y)\\
				\leqslant&\frac{C}{\ve^2}\int_{B_{8\ve}(x)}\frac{1}{|x-y|}\ud V_{\S^2}(y)\leqslant \frac{C}{\ve}.
			\end{align*}
			
			Notice that $\nabla u\in H^{1/2}(\S^2)$ for any $u\in H^{3/2}(\S^2)$. It suffices to show
			  that given $u\in H^{1/2}(\S^2)$, there holds
	\begin{align}\label{approximate lem claim}
		\lim_{\ve\to 0}\underbrace{\int_{\S^2} \int_{\S^2} \frac{\eta_{\ve}(x)\eta_{\ve}(y)(u(x)-u(y))^2}{|x-y|^3}  \ud V_{\S^2}(y)  \ud V_{\S^2}(x)}_{I^{1}_{\e}}=0.
	\end{align}
	
	For $u\in C^{\infty}(\S^2)$ we have
	\begin{align}\label{approximate lem formula 3}
		I^{1}_{\ve}\leqslant C\|u\|^2_{C^{1}(\S^2)}\int_{B_{2\ve}(S)}\int_{B_{2\ve}(S)}\frac{1}{|x-y|}\ud V_{\S^2}(y)  \ud V_{\S^2}(x)\leqslant C\|u\|^2_{C^{1}(\S^2)} \ve^3.
	\end{align}
	Since $C^{\infty}(\S^2)$ is dense in $H^{1/2}(\S^2)$, given $u \in H^{1/2}(\S^2)$ there exists a sequence $\{u_n\}\subset C^\infty(\S^2)$ such that $\|u_n-u\|_{H^{1/2}(\S^2)}\to 0$. If we let $v_n=u-u_n$, then
	\begin{align}\label{approximate lem formula 4}
		\lim_{n\to +\infty}\int_{\S^2} \int_{\S^2} \frac{(v_n(x)-v_n(y))^2}{|x-y|^3}  \ud V_{\S^2}(y)  \ud V_{\S^2}(x)=0.
	\end{align}
	
	By \eqref{approximate lem formula 3} we have
	\begin{align*}
		&\int_{\S^2} \int_{\S^2} \frac{\eta_{\ve}(x)\eta_{\ve}(y)(u(x)-u(y))^2}{|x-y|^3}  \ud V_{\S^2}(y)  \ud V_{\S^2}(x)\\
		\leqslant &2	\int_{\S^2}\int_{\S^2} \frac{\eta_{\ve}(x)\eta_{\ve}(y)(v_n(x)-v_n(y))^2}{|x-y|^3}  \ud V_{\S^2}(y)  \ud V_{\S^2}(x)\\
		&+2	\int_{\S^2}\int_{\S^2} \frac{\eta_{\ve}(x)\eta_{\ve}(y)(u_n(x)-u_n(y))^2}{|x-y|^3}  \ud V_{\S^2}(y)  \ud V_{\S^2}(x)\\
		\leqslant& 2\int_{\S^2}\int_{\S^2} \frac{(v_n(x)-v_n(y))^2}{|x-y|^3}  \ud V_{\S^2}(y)  \ud V_{\S^2}(x)+C \|u_n\|_{C^{1}(\S^2)}^2 \ve^3.
	\end{align*}
Thus, we obtain
	\begin{align}\label{est:limsup}
		\limsup_{\ve\to 0}I^{1}_{\ve}\leqslant 	\int_{\S^2} \int_{\S^2} \frac{(v_n(x)-v_n(y))^2}{|x-y|^3}  \ud V_{\S^2}(y)  \ud V_{\S^2}(x).	
		\end{align}
	Next let $n \to \infty$, the claim follows by \eqref{est:limsup} and  \eqref{approximate lem formula 4}. 
	
	\medskip
	
	Therefore, going back to \eqref{approximate lem formula 2}, by \eqref{approximate lem claim} we have
	\begin{align*}
		I_{\ve}\leqslant C\|\nabla u\|^2_{L^4(B_{2\ve}(S))}+o_{\ve}(1)\to 0 \qquad \mathrm{~~as~~} \ve\to 0.
	\end{align*}

	For the second term, similarly we have
	\begin{align*}
		II_{\ve}:=&\int_{\S^2}\partial_i\eta_{\ve}(x)u(x)\int_{\S^2} \frac{\partial_i\eta_{\ve}(x)u(x)-\partial_i\eta_{\ve}(y)u(y)}{|x-y|^3}  \ud V_{\S^2}(y)  \ud V_{\S^2}(x)\\
	=&\int_{\S^2}\partial_i\eta_{\ve}(x)u^2(x)\int_{\S^2} \frac{\partial_i\eta_{\ve}(x)-\partial_i\eta_{\ve}(y)}{|x-y|^3}  \ud V_{\S^2}(y)  \ud V_{\S^2}(x)\nonumber\\
	&+\frac{1}{2}\int_{\S^2} \int_{\S^2} \frac{\partial_i\eta_{\ve}(x)\partial_i\eta_{\ve}(y)(u(x)-u(y))^2}{|x-y|^3}  \ud V_{\S^2}(y)  \ud V_{\S^2}(x)\\
	:=&II^{1}_{\ve}+II^{2}_{\ve}.
	\end{align*}
	For  $II^1_{\ve}$, by \eqref{lim:intermediate_norm} we estimate
	\begin{align*}
		|II^{1}_{\ve}|\leqslant& \frac{C}{\ve}\|\nabla^3\eta_{\ve}\|_{L^{\infty}(\S^2)}\int_{B_{2\ve}\backslash B_{\ve}(S)}u^2(x)\int_{B_{4\ve}(S)}\frac{1}{|x-y|}\ud V_{\S^2}(y)\ud V_{\S^2}(x)\\
		\leqslant &\frac{C[u]_{C^{1/2}(B_{2\ve}(S))}^2}{\ve^3}\int_{B_{2\ve}\backslash B_{\ve}(S)}d_{\S^2}(S,x)\ud V_{\S^2}(x)\\
		\leqslant& C[u]_{C^{1/2}(B_{2\ve}(S))}^2\to 0.
	\end{align*}
	For  $II^2_{\ve}$, we have
	\begin{align*}
		II^2_{\ve}\leqslant& \frac{C}{\ve^2}\int_{B_{2\ve}\backslash B_{\ve}(S)}\int_{B_{2\ve}\backslash B_{\ve}(S)} \frac{(u(x)-u(y))^2}{|x-y|^3}  \ud V_{\S^2}(y)  \ud V_{\S^2}(x)\\
		\leqslant& \frac{C}{\ve^2}\int_{\S^2} \int_{\S^2} \frac{(\eta_{2\ve}(x)u(x)-\eta_{2\ve}(y)u(y))^2}{|x-y|^3}  \ud V_{\S^2}(y)  \ud V_{\S^2}(x)\\
		\leqslant& \frac{C}{\ve^2}\|\eta_{2\ve}u\|^2_{H^{1/2}(\S^2)}.
	\end{align*}
	 The interpolation inequality \eqref{ineq:interpolation} gives
	\begin{align*}
		\|\eta_{2\ve}u\|^2_{H^{1/2}(\S^2)}\leqslant C\|\eta_{2\ve}u\|_{L^2(\S^2)} \|\eta_{2\ve}u\|_{H^1(\S^2)}.
		\end{align*}
		Observe that
		\begin{align*}
	 \|\eta_{2\ve}u\|_{L^2(\S^2)}\leqslant& C[u]_{C^{1/2}(B_{2\ve}(S))}\left(\int_{B_{2\ve}\backslash B_{\ve}(S)}d_{\S^2}(S,x)\ud V_{\S^2}(x)\right)^{1/2}\\
			\leqslant& C[u]_{C^{1/2}(B_{2\ve}(S))}\ve^{3/2}
	\end{align*}
and
	\begin{align*}
	\left(\int_{\S^2} |\nabla (\eta_{\ve} u)|^2  \ud V_{\S^2}\right)^{1/2}\leqslant& \frac{C}{\ve}\left(\int_{B_{2\ve}(S)}u^{2}\ud V_{\S^2}\right)^{1/2}+C\left(\int_{B_{2\ve}(S)}|\nabla u|^{2}\ud V_{\S^2}\right)^{1/2}\\
	\leqslant&C[u]_{C^{1/2}(B_{2\ve}(S))}\ve^{1/2}+C\left(\int_{B_{2\ve}(S)}|\nabla u|^{4}\ud V_{\S^2}\right)^{1/4}\ve^{1/2}.
	\end{align*}
	Thus, we obtain
	\begin{align*}
		\|\eta_{2\ve}u\|_{H^{1/2}(\S^2)}\leqslant C\|u\|_{W^{1,4}(B_{4\ve}(S))}\ve,
	\end{align*}
	whence,
	\begin{align*}
		II^2_{\e}\leqslant C\|u\|^2_{W^{1,4}(B_{4\e}(S))}\to 0.
	\end{align*}

Finally, combining the estimates above we deduce that
	\begin{align*}
		\|\partial_i(\eta_{\ve}u)\|^2_{H^{1/2}(\S^2)}\leqslant I_{\ve}+II_{\ve}+\int_{\S^2} |\nabla (\eta_{\ve} u)|^2  \ud V_{\S^2}\to 0
	\end{align*}
	and thus $\eta_{\ve} u\to 0$ in $H^{3/2}(\S^2)$.
		\end{pf}

		The following result is straightforward, so we omit the proof.
		\begin{lem}\label{V3 lem 3}
			Suppose $u\in H^{3/2}(\S^2)$, there exist positive constants $C_1$ and $C$ such that \begin{align*}
				\int_{\S^2}uP_3u   \ud V_{\S^2}+C_1\int_{\S^2} u^2  \ud V_{\S^2}\geqslant C\|u\|^2_{H^{3/2}(\S^2)}.
			\end{align*}
		\end{lem}

		Very recently, it has been known in the first author and Shi \cite[Theorem 1 (3)]{Chen-Shi}  that the representation formula of Green function for $P_3$ on $\S^2$ is 
		\begin{equation}\label{Green_fcn}
		G_{x_0}(\cdot)=-\frac{1}{2\pi} |\cdot - x_0| \qquad \mathrm{~~for~~some~~} x_0 \in \S^2
		\end{equation}
		satisfying $P_3 G_{x_0}=\delta_{x_0}$ in the distribution sense, where $ |\cdot - x_0|$ means the Eclidean distance from $x_0$ in $\R^3$.

	\begin{prop}\label{Prop: nonnegative_energy}
		Suppose $u\in H^{3/2}(\S^2)$ and $u(x_0)=0$ for some $x_0\in \S^2$, then $E[u]\geqslant 0$. Assume additionally that $u\geqslant 0$, then  $E[u]=0$ if and only if $u=CG_{x_0}$ for some $C \in \R_-$.
	\end{prop}	
	\begin{pf}
		Without loss of generality, we assume $x_0=S$. By Lemma \ref{Approximate lem}, there exists a sequence $\{u_n\}\subset C^{\infty}(\S^2)$ such that $u_n$ vanishes near $S$ and $\|u_n-u\|_{H^{3/2}(\S^2)}\to 0$ as $n\to+\infty$. Set $v(x)=u\circ I(x)\sqrt{\frac{1+|x|^2}{2}}$ and $v_n(x)=u_n\circ I(x)\sqrt{\frac{1+|x|^2}{2}}\in C^{\infty}_{c}(\R^2)$. Then by Proposition \ref{GJMS:extrinsic_intrinsic} and Theorem \ref{Dim 3 integral equ thm} combined with the estimate $|E[u]|\leqslant C\|u\|^2_{H^{3/2}(\S^2)}$, we obtain
		\begin{align*}
			\int_{\R^2}\left(\left(-\Delta\right)^{3/4}(v_n-v_m)\right)^2 \ud x=E[u_n-u_m]\to 0\qquad \mathrm{as}\qquad n,m\to+\infty.
		\end{align*}
		This indicates that $\left(-\Delta\right)^{3/4} v_n$ is a Cauchy sequence in $L^2(\R^2)$. Then there exists  $f\in L^2(\R^2)$ such that 
		\begin{align*}
			\left\|\left(-\Delta\right)^{3/4} v_n-f\right\|_{L^2(\R^2)}\to 0 \qquad \mathrm{as}\qquad n\to+\infty.
		\end{align*}
		For any $\varphi\in C^{\infty}_{c}(\R^2)$ we have
		\begin{align*}
			\int_{\R^2}\left(-\Delta\right)^{3/4} v_n\varphi \ud x=	\int_{\R^2} v_n\left(-\Delta\right)^{3/4}\varphi\ud x\to \int_{\R^2} v\left(-\Delta\right)^{3/4}\varphi \ud x\qquad \mathrm{as}\qquad n\to+\infty.
		\end{align*}
		On the other hand, 
	\begin{align*}
		\int_{\R^2}\left(-\Delta\right)^{3/4} v_n\varphi \ud x \to \int_{\R^2}f\varphi \ud x\qquad \mathrm{~~as~~}\quad n\to+\infty.
	\end{align*}
	This means $\left(-\Delta\right)^{3/4}v=f$ in the distribution sense. 
	
	Observe that
	\begin{align*}
		\left|E[u_n]-\int_{\R^2}f^2 \ud x\right|=&\left|\int_{\R^2}\left(\left(-\Delta\right)^{3/4} v_n-f\right)\left(\left(-\Delta\right)^{3/4} v_n+f\right)\ud x\right|\\
		\leqslant&\left\|\left(-\Delta\right)^{3/4} v_n-f\right\|_{L^2(\R^2)}\left\|\left(-\Delta\right)^{3/4} v_n+f\right\|_{L^2(\R^2)}\\
		\leqslant& C\left\|\left(-\Delta\right)^{3/4} v_n-f\right\|_{L^2(\R^2)}\to 0.
	\end{align*}
	By  \eqref{Energy inequlity} and  the fact that up to a subsequence, $u_n\to u$ in $C^{\theta}(\S^2), ~\forall~ \theta \in (0,1/2)$, we have
	 \begin{align*}
		\left|E[u]-E[u_n]\right|\leqslant&\left|\int_{\S^2} (u-u_n)P_3u_n  \ud V_{\S^2}+\int_{\S^2}u P_3 (u-u_n)  \ud V_{\S^2}\right|\\
		\leqslant& \left(E[u-u_n]+\frac{3}{8}|\overline{u-u_n}|^2\right)^{1/2}\left(E[u]+\frac{3}{8}\bar{u}^2\right)^{1/2}+\frac{3}{8}|\overline{u-u_n}|(|\bar{u}|+|\overline{u_n}|)\\
		&+ \left(E[u-u_n]+\frac{3}{8}|\overline{u-u_n}|^2\right)^{1/2}\left(E[u_n]+\frac{3}{8}|\overline{u_n}|^2\right)^{1/2}\\
		\to &0.
	\end{align*}
	
	Hence, we obtain
	\begin{align*}
		E[u]=\int_{\R^2}f^2 \ud x=\int_{\R^2}\left(\left(-\Delta\right)^{3/4}v\right)^2 \ud x\geqslant 0.
	\end{align*}

	Now assume  $u\geqslant 0$ and $E[u]=0$,  then
	\begin{align*}
		\left(-\Delta\right)^{3/4}v=0\qquad\mathrm{in}\qquad \R^2
	\end{align*}
	in the distribution sense.
	Notice that $u(S)=0$ and $\|u\|_{C^{1/2}(\S^2)}\leqslant C$, then for $|x|\gg 1$, 
	\begin{align*}
	0\leqslant u\circ I(x)\leqslant C|x|^{-1/2}
	\end{align*}
	and thus $0\leqslant v(x)\leqslant C|x|^{1/2}$. 
 By the Liouville theorem in \cite{Zhou&Chen&Cui&Yuan}, we know $v\equiv C$, that is,
	\begin{align*}
		u\circ I(x)=C\sqrt{\frac{2}{1+|x|^2}}.
	\end{align*}
	This together with \eqref{Green_fcn} implies
	$$u=-2\sqrt{2}\pi C G_{x_0}.$$
	
	The opposite direction follows from \eqref{Green_fcn}.
	\end{pf}	
	
	The following is a simple observation: If there exists  $0\leqslant u \in H^{3/2}(\S^2)$ vanishing somewhere such that $E[u]<0$, then for any $\ve>0$, we have $E[u+\ve] \to E[u]$,  and a contradiction argument together with Lemmas \ref{V3 lem 1} and \ref{V3 lem 3} yields $\|(u+\ve)^{-1}\|_{L^4(\S^2)} \to \infty$  as $\ve \to 0$, thus $Y_3^+(\S^2)=-\infty$. In this sense, Proposition \ref{Prop: nonnegative_energy} provides a necessary condition of $Y_3^+(\S^2)$ to be finite.
	
	\medskip
	
	\noindent \textbf{Proof of Theorem \ref{thm:Sobolev_ineq}.} 
	By conformal covariance we can find a sequence of positive functions $\{u_n\} \subset H^{3/2}(\S^2)$   such that 
		\begin{align}\label{Variation lem formula a}
		\max_{\S^2}u_n=1 \qquad \mathrm{~~and~~}\qquad\qquad E[u_n]\|u_n^{-1}\|^2_{L^{4}(\S^2)}\to Y_3^+(\S^2).
	\end{align} 
	We point out that at this stage $Y_3^+(\S^2)$ might be $-\infty$.
	Then $\|u_n^{-1}\|_{L^{4}(\S^2)}\geqslant |\S^2|^{1/4}$ and 
	\begin{align*}
		E[u_n]\|u_n^{-1}\|^2_{L^{4}(\S^2)}\leqslant C.
	\end{align*}
	By Lemma \ref{V3 lem 3} we have
	\begin{align*}
		\|u_n\|_{H^{3/2}(\S^2)}^2-C\|u_n\|_{L^{2}(\S^2)}^2\leqslant E[u_n]\leqslant C.
	\end{align*}
	This implies  $\|u_n\|_{H^{3/2}(\S^2)}\leqslant C$. Then up to a subsequence we have
	\begin{align*}
		u_n \rightharpoonup u_{\infty} \quad \mathrm{in}\quad H^{3/2}(\S^2) \qquad \mathrm{and}\qquad u_n \rightarrow u_{\infty} \quad \mathrm{~~in~~}\quad C^{\theta}(\S^2),
	\end{align*}
	where $\theta\in (0,1/2)$. 
	
	Now two possibilities of $u_\infty \geqslant 0$ occur.
	
	\medskip
	 \noindent\underline{\emph{Case 1.}}  $u_{\infty}>0$  on $\S^2$. 
	\medskip
	
		We have  $u_n^{-1}\to u_{\infty}^{-1}$ uniformly on $\S^2$, and $\|u_n^{-1}\|_{L^{4}(\S^2)}\to \|u_{\infty}^{-1}\|_{L^4(\S^2)}$. Then 
	\begin{align*}
	E[u_{\infty}]\|u_{\infty}^{-1}\|^2_{L^{4}(\S^2)}\leqslant \limsup_{n\to +\infty}E[u_n]\|u_n^{-1}\|^2_{L^{4}(\S^2)}=Y_3^+(\S^2).
	\end{align*}
	This implies
	\begin{align*}
	Y_3^+(\S^2)=E[u_{\infty}]\|u_{\infty}^{-1}\|^2_{L^{4}(\S^2)}.
	\end{align*}
	
		\medskip
	 \noindent\underline{\emph{Case 2.}} $u_{\infty}(x_0)=0$ for some $x_0 \in \S^2$. 
		\medskip
		
	Without loss of generality, we assume $x_0=S$. It follows from Lemma \ref{V3 lem 1} that $\|u_{\infty}^{-1}\|_{L^{4}(\S^2)}=+\infty$, then
	\begin{align*}
		+\infty=\|u_{\infty}^{-1}\|_{L^{4}(\S^2)}\leqslant \liminf_{n\to +\infty}\|u_n^{-1}\|_{L^{4}(\S^2)}.
	\end{align*}
	This implies
	\begin{align*}
		E[u_{\infty}]\leqslant \limsup_{n\to +\infty}E[u_n]\leqslant 0.
	\end{align*}
	This together with Proposition \ref{Prop: nonnegative_energy} yields $E[u_{\infty}]=0$  and $u_{\infty}=-\pi G_{x_0}$, that is, $u_\infty \circ I(y)=(1+|y|^2)^{-1}$.  On the other hand, we can find $x_n \in \S^2$ such that $u_n(x_n)=\min_{\S^2}u_n:= \lambda_n$, and $\lim_{n\to +\infty}x_n=S$.
	
	Denote by $I_{x_n}:\R^2\to \S^2\backslash\{-x_n\} $ and $\delta_\lambda(y)=\lambda y$ for $y \in \R^2$, where $\lambda \in \R_+$.  We consider a conformal transformation on $\S^2$:
	$$\varphi_{\lambda_n}=I_{x_n}\circ\delta_{\lambda_n}\circ I^{-1}_{x_n}.$$
	We define
	\begin{align}\label{Variation lem formula b}
		v_n\circ I_{x_n}(y)=\left(\frac{1+\lambda^2_n|y|^2}{\lambda_n(1+|y|^2)}\right)^{1/2}u_n\circ I_{x_n}(\lambda_ny)
	\end{align}
	such that
	$v_n^{-4}g_{\S^2}=\varphi_{\lambda_n}^\ast\left(u_n^{-4}g_{\S^2}\right)$. Then, it follows from conformal covariance of $P_3$ that $v_n$ is another minimizing sequence of $Y_3^+(\S^2)$.
	
		By definition \eqref{Variation lem formula b} of $v_n$ we have
	\begin{align}\label{Variation lem formula c}
		v_n(x_n)=v_n\circ I_{x_n}(0)=\sqrt{\lambda_n}\qquad \mathrm{and}\qquad \lim_{n\to +\infty} \frac{v_n(N)}{\sqrt{\lambda_n}}=1.
	\end{align}
	Let $\nu_n=\max_{\S^2} v_n$ and $\tilde{v}_n:=\frac{v_n}{\nu_n}$, then up to a subsequence, the same argument as above yields 
	\begin{align*}
		\tilde{v}_n \rightharpoonup \tilde{v}_{\infty} \quad \mathrm{in}\quad H^{3/2}(\S^2) \qquad \mathrm{and}\qquad \tilde{v}_n \rightarrow \tilde{v}_{\infty} \quad \mathrm{in}\quad C^{\theta}(\S^2),
	\end{align*}
	 	where $\theta\in (0,1/2)$. 	By \eqref{Variation lem formula c} and Proposition \ref{Prop: nonnegative_energy}, we must have $\tilde{v}_{\infty}(N)=\tilde{v}_{\infty}(S)>0$, whence $\tilde{v}_n(N)=\frac{\sqrt{\lambda_n}}{\nu_n}\to \tilde{v}_{\infty}(N)$.   On the other hand, for all $y \in I^{-1}(\S^2\backslash\{S,N\})$ we have
	\begin{align*}
		\tilde{v}_n\circ I_{x_n}(y)=&\frac{v_n\circ I_{x_n}(y)}{\nu_n}=\left(\frac{1+\lambda^2_n|y|^2}{\lambda_n(1+|y|^2)}\right)^{1/2}\frac{u_n\circ I_{x_n}(\lambda_ny)}{\nu_n}\\
		\geqslant&\frac{\lambda_n}{\nu_n}\left(\frac{1+\lambda^2_n|y|^2}{\lambda_n(1+|y|^2)}\right)^{1/2}\to \frac{\tilde{v}_{\infty}(N)}{(1+|y|^2)^{1/2}}>0.
	\end{align*}
	Hence, we arrive at $\tilde{v}_{\infty}>0$ on $\S^2$ and also 
	\begin{align*}
		Y_3^+(\S^2)=E[\tilde{v}_{\infty}]\|\tilde{v}_{\infty}^{-1}\|^2_{L^{4}(\S^2)}.
	\end{align*}

Therefore, the minimizer, a positive function $u \in H^{3/2}(\S^2)$, is achieved and satisfies, modulo  a positive constant,
\begin{align*}
	P_3u=-\frac{3}{8}u^{-5}.
\end{align*}
Set $I:\R^2\to \S^2\backslash\{S\}$ and 
\begin{align}
v(y)=u\circ I(y)\sqrt{\frac{1+|y|^2}{2}}.
\end{align}
By Proposition \ref{GJMS:extrinsic_intrinsic}  and  Theorem \ref{Dim 3 integral equ thm}, we know
\begin{align*}
	v(y)=\frac{3}{16\pi}\int_{\R^2}|y-z|v^{-5}(z)\ud z.
\end{align*}
		By the classification theorem for conformally invariant integral equations in Yan Yan Li \cite[Theorem 1.5]{Li}, we obtain
		\begin{align*}
			v(y)=\sqrt{\frac{\ve^2+|y-y_0|^2}{2\ve}} \qquad\mathrm{~~for~~some~~}\quad y_0\in \R^2,\quad\ve\in \R_{+}.
		\end{align*}
		Hence, the representation of solution on $\S^2$ follows by an appropriate stereographic projection.
		\hfill $\Box$
		
		\section{Extremal functions on balls for Sobolev trace inequalities}\label{Sect6}
		
		Recently, Ndiaye and L. Sun \cite{Ndiaye-Sun} invoked an integral equation method to study extremal functions on balls for Ache-Chang's Sobolev trace inequalities. Our purpose is to make a geometric interpretation of these biharmonic extremal functions on balls, as well as the ones for ours. Our approach is of geometric favor, based on the explicit formula of extremal function on spheres due to Theorem \ref{thm:Sobolev_ineq} and Theorem A, and thus is more straightforward. We emphasize the importance of boundary defining function on $\Bn$ .\footnote{ This idea originates from ` \emph{a toy model} ' , which is set up for examples of optimal constants  by the first author, Wei and Wu \cite[Section 2]{Chen-Wei-Wu}.}

		For $n\geqslant 3$, fix  $a \in \Bn$ and set $\bar a=|a|^{-2}a$ whenever $a \neq 0$.  We define by a conformal transformation on $\overline{\Bn}$ by
\begin{align}\label{conf_trans_ball_Hua}
y=\psi_a(x)=&\frac{x-a-|x|^2 a+2(a \cdot x)a-|a|^2 x}{1-2 a \cdot x+|a|^2|x|^2}\\
=&\frac{x-a-|x-\bar a|^2 a-2(\bar a \cdot x) a+|\bar a|^2 a+2(a \cdot x)a-|a|^2 x}{|a|^2|x-\bar a|^2} \qquad \mathrm{~~for~~every~~} a \neq 0\no\\
=&-\bar a+\frac{(1-|a|^2)(x-2(\bar a \cdot x)a+\bar a)}{|a|^2|x-\bar a|^2}\no\\
=&-\bar a+\frac{(|\bar a|^2-1)(x-2(\bar a \cdot x)a+\bar a)}{|x-\bar a|^2} \rightarrow x \quad \mathrm{~~as~~} a \to 0,\no
\end{align}
with the property that
$$\psi_a^\ast (|\ud y|^2)=\left(\frac{1-|a|^2}{|x|^2 |a|^2-2a \cdot x+1}\right)^2 |\ud x|^2$$
and $\psi_a (\Bn)=\Bn, \psi_a(\pa \Bn)=\pa \Bn$. Moreover we have
$$\frac{|\ud y|^2}{(1-|y|^2)^2}=\frac{|\ud x|^2}{(1-|x|^2)^2}.$$
In other words, $\psi_a$ is also an isometry on the Poincar\'e ball $(\Bn,(\frac{2}{1-|x|^2})^2 |\ud x|^2)$. See L.-K. Hua \cite[Chapter 1]{Hua}. The following is a geometric interpretation of $\psi_a$: If we denote by
 $$\varphi_a (x)=\frac{(x-\bar a)(|\bar a|^2-1)}{|x-\bar a|^2}+\bar a, \quad a \neq 0$$
 an inversion with respect to a sphere with radius $\sqrt{|\bar a|^2-1}$ and centered at $\bar a$, then
 $$\psi_a(x)=\varphi_{- a}(x-2(\bar a \cdot x)a).$$
 
 Suppose $F$ is an inversion with respect to the sphere $\pa B_{\sqrt 2}(-\mathbf{e}_n)$, which conformally maps from $(\mathbb{R}^n_+,|\ud z|^2)$ to $(\Bn, |\ud x|^2)$, explicitly,
$$x=F(z)=-\mathbf{e}_n+\frac{2(z+\mathbf{e}_n)}{|y+\mathbf{e}_n|^2}=\frac{(2z',1-|z|^2)}{(1+z_n)^2+|z'|^2}$$
for $z=(z',z_n) \in\mathbb{R}^n_+$.

Now we define
$$W_{(z_0',\ve)}(z)=\log\frac{2\ve}{(\ve+z_n)^2+|z'-z'_0|^2}$$
and
$$\hat U_a(x)=\log \frac{1-|a|^2}{|a|^2 |x|^2-2a \cdot x+1}\qquad \mathrm{~~such~~that~~}\quad e^{2\hat U_a(x)}|\ud x|^2=\psi_a^\ast(|\ud y|^2).$$
One can verify that
\begin{equation}\label{bubbles_ball_half_space}
F^\ast(e^{2\hat U_a(x)}|\ud x|^2)=e^{2W_{(z_0',\ve)}(z)} |\ud z|^2 \qquad \mathrm{for~~some~~} \ve \in \R_+, z_0' \in \R^{n-1}, 
\end{equation}
via the change of parameters
\begin{equation}\label{formula_parameters}
\ve=\frac{1-|a|^2}{|a|^2+2a_n+1}, \quad z_0'=\frac{2a'}{|a|^2+2a_n+1} \quad \Longleftrightarrow \quad (z_0',\ve)=F^{-1}(a)=F(a).
\end{equation}
For scalar flat conformal metrics with constant boundary mean curvature problem, the above formulae \eqref{bubbles_ball_half_space} and \eqref{formula_parameters} represent the relationship between standard bubbles in $\Bn$ and $\R^n_+$.
				
		We  first restrict consideration to $n\geqslant 3$ and $n \neq 4$. A similar argument as  Proposition \ref{prop:Kazdan-Warner} shows that modulo a positive constant, a smooth positive minimizer of \eqref{ineq:Ache-Chang_n>4} solves
		\begin{align}\label{EL-eqn:Ache-Chang_dim_neq_four}
	\begin{cases}
	\displaystyle \Delta^2 U=0  &\mathrm{~~in~~}\qquad \Bn,\\
	\displaystyle 	\frac{\partial U}{\partial r}=-\frac{n-4}{2}u &\mathrm{~~on~~}\qquad \S^{n-1},\\
	\displaystyle \mathscr{B}_3^3 U= \frac{n(n-2)(n-4)}{4} u^{\frac{n+2}{n-4}} \qquad &\mathrm{~~on~~}\qquad \S^{n-1},
	\end{cases}
	\end{align}
	where $\mathscr{B}_3^3$ is given in \eqref{bdry operator_dim_n}.
		
		Let
$$\hat U_a(x)=\Big(\frac{1-|a|^2}{|a|^2 |x|^2-2a \cdot x+1}\Big)^{\frac{n-4}{2}}$$
such that $\psi_a^\ast(|\ud x|^2)=\hat U_a(x)^{\frac{4}{n-4}}|\ud x|^2$.  We introduce
$$U_a(x)=\hat U_a(x)+\frac{1-|x|^2}{2}h(x) \qquad\mathrm{~~for~~some~~} h\in C^\infty(\overline {\Bn}).$$
Here $g$ is chosen to fulfill the Neumann boundary condition in \eqref{cond:Neumann_dim_three} and \eqref{cond:Neumann_dim_above_three}
\begin{align*}
0=\frac{\pa U_a}{\pa r}+\frac{n-4}{2}U_a=&\frac{\pa \hat U_a}{\pa r}+\frac{n-4}{2}\hat U_a-h(x) \\
=&\frac{n-4}{2}\hat U_a(x) \frac{1-|a|^2 |x|^2}{|a|^2 |x|^2-2a \cdot x+1}-h(x)\qquad \qquad \mathrm{on~~}\quad \S^{n-1}.
\end{align*}
Besides the above condition, meanwhile requiring that $(1-|x|^2)h(x)$ is biharmonic in $\overline{\Bn}$, we may choose
$$h(x)=\frac{n-4}{2}\Big(\frac{1-|a|^2}{|a|^2 |x|^2-2a \cdot x+1}\Big)^{\frac{n-2}{2}}.$$
Thus, for each $a \in \Bn$,
\begin{equation}\label{extremal_fcns_on_balls}
U_a(x)=\Big(\frac{1-|a|^2}{|a|^2 |x|^2-2a \cdot x+1}\Big)^{\frac{n-4}{2}}+\frac{n-4}{4}(1-|x|^2)\Big(\frac{1-|a|^2}{|a|^2 |x|^2-2a \cdot x+1}\Big)^{\frac{n-2}{2}}
\end{equation}
forms a family of smooth solutions on $\overline{\Bn}$ to PDE \eqref{EL-eqn:Ache-Chang_dim_neq_four}.

\medskip

The same trick also works for  extremal functions for the  Sobolev trace inequality  by Ache-Chang \cite{Ache-Chang} on $\B^4$, which is stated as follows.
\begin{thm B}
Given  $u\in C^{\infty}(\mathbb{S}^{3})$, let $U$ be a smooth extension of $u$ to  $\mathbb{B}^4$ satisfying
\begin{align*}
	\frac{\pa U}{\partial r}=0 \qquad \mathrm{~~on~~} \S^3.
\end{align*}
 Then with $\bar u:=\int_{\S^3}u \ud \mu_{\S^3}/(2\pi^2)$, there holds
\begin{align}\label{ineq:Ache-Chang_n=4}
	\log \left(\frac{1}{2 \pi^{2}} \int_{\S^{3}} e^{3(u-\bar{u})} \ud V_{\S^3}\right)
	 \leqslant \frac{3}{16 \pi^{2}}\left[ \int_{\mathbb{B}^{4}}\left(\Delta U\right)^{2} \ud x+2 \int_{\S^{3}}|\nabla u|_{\S^3}^{2} \ud V_{\S^3}\right],
\end{align}
with equality if and only if $U$ is a biharmonic extension of some function $u_{z_0}(x)=-\log|1-z_0\cdot x|+C$ on $\S^{3}$, and satisfies zero Neumann boundary condition, where $z_0 \in \mathbb{B}^4, C \in \R$. 
\end{thm B}

It is not hard to see that modulo a constant, a smooth minimizer of the inequality \eqref{ineq:Ache-Chang_n=4} solves
	\begin{align}\label{EL-eqn:Ache-Chang_dim_four}
	\begin{cases}
	\displaystyle \Delta^2 U=0  &\mathrm{~~in~~}\qquad \B^4,\\
	\displaystyle 	\frac{\partial U}{\partial r}=0 &\mathrm{~~on~~}\qquad \S^3,\\
	\displaystyle \mathscr{B}_3^3 U+2=2e^{3u} \qquad &\mathrm{~~on~~}\qquad \S^3,
	\end{cases}
	\end{align}
	where $\mathscr{B}_3^3$ is given in \eqref{bdry operator_dim_n} for $n=4$.

To fulfill the vanishing Neumann boundary condition, we introduce
$$U_a(x):=\hat U_a(x)+\frac{1-|x|^2}{2}\hat h(x) \qquad\mathrm{~~for~~some~~} \hat h\in C^\infty(\overline {\B^4})$$
such that
\begin{align*}
0=\frac{\pa U_a}{\pa r}=\frac{\pa \hat U_a}{\pa r}-\hat h(x) \qquad \mathrm{on~~}\quad \S^3.
\end{align*}
Besides the above condition, ensuring that $(1-|x|^2)\hat h(x)$ is biharmonic in $\overline{\B^4}$, forces us to take
$$\hat h(x)=\frac{1-|a|^2}{|a|^2 |x|^2-2a \cdot x+1}-1,$$
whence,
$$U_a(x)=\log \frac{1-|a|^2}{|a|^2 |x|^2-2a \cdot x+1}+\frac{1-|x|^2}{2}\big(\frac{1-|a|^2}{|a|^2 |x|^2-2a \cdot x+1}-1\big), \quad a \in \B^4$$
forms the set of smooth solutions in $\overline{\B^4}$ to PDE \eqref{EL-eqn:Ache-Chang_dim_four}.

Combining \eqref{bubbles_ball_half_space} and \eqref{formula_parameters} we deduce that
$$F^\ast(e^{2U_a(x)}|\ud x|^2)=e^{2 V_{(z_0',\ve)}(z)} |\ud z|^2.$$

\medskip

\noindent\textbf{Completion of the proof of Theorem \ref{Thm:Sobolev_trace_ineq_three_balls}.} As indicated in Section \ref{Sect3}, the combination of Theorem \ref{thm:Sobolev_ineq} and the representation formula \eqref{extremal_fcns_on_balls} of extremal functions for $n=3$ is enough to finish the proof of Theorem \ref{Thm:Sobolev_trace_ineq_three_balls}. \hfill $\Box$

\medskip

As an application, we include a sharp Sobolev trace inequality on $\B^3$ without constraints.
		
			\begin{cor}
		 Suppose $U\in C^\infty(\overline{\B^3})$ with $U\big|_{\S^2}=u>0$. Let 
		 \begin{align}\label{fcn:amended}
	\hat{U}=\frac{|x|^2+3}{4}U+\frac{1-|x|^2}{2}x\cdot\nabla U+V,
	\end{align}
	for all $V\in C^{\infty}(\B^3)$ satisfying
	\begin{align*}
		V=0 \qquad \mathrm{and}\qquad \frac{\partial V}{\partial r}=0  \qquad \mathrm{on}\qquad \S^2.
	\end{align*}
		 Then
		\begin{align*}
		-\frac{3}{4}|\mathbb{S}^{2}|^{\frac{3}{2}} \left(\int_{\mathbb{S}^{2}}|u|^{-4}\ud V_{\mathbb{S}^{2}}\right)^{-\frac{1}{2}}
		\leqslant \int_{\mathbb{B}^3} \left(\Delta\hat{U}\right)^2 \ud x + 2\int_{\mathbb{S}^{2}}|\nabla u|_{\S^{2}}^2 \ud V_{\S^{2}} -\frac{3}{2}\int_{\mathbb{S}^{2}} u^2 \ud V_{\S^{2}},
		\end{align*}
		Moreover, equality holds if and only if $\hat U$ is biharmonic on $\overline{\B^3}$ as in Theorem \ref{Thm:Sobolev_trace_ineq_three_balls}.
	\end{cor}
	
	\begin{pf}
	A direct computation yields
	\begin{align*}
		\hat{U}=u \qquad\qquad \mathrm{and}\qquad\qquad \frac{\partial \hat{U}}{\partial r}=\frac{u}{2} \qquad \mathrm{~~on~~}\qquad \S^2.
	\end{align*}
	Then, this is a direct consequence of Theorem \ref{Thm:Sobolev_trace_ineq_three_balls}.
	\end{pf}
		\begin{rem}
		The choice of $V$ in \eqref{fcn:amended} could be wide.  For instance,  $V$ is every function in $C^{\infty}_{c}(\B^3)$,  or  $V=(1-|x|^2)^{2k}W(x)$ for some $k\in \mathbb{Z}_+$ and $W\in C^{\infty}(\overline{\B^3})$, etc. In particular, taking $V=0$ we obtain a fourth order sharp Sobolev trace inequality without constraints on $\B^3$.
		\end{rem}

		\bibliographystyle{unsrt}

	\end{document}